\documentclass[12pt]{amsart}
\usepackage{amsmath,amsfonts,latexsym,graphicx,amssymb,url,geometry}
\usepackage{hyperref,pdfsync}
\usepackage{amsmath,txfonts,pifont,bbding,pxfonts,manfnt}
\usepackage[active]{srcltx}
\usepackage{wasysym,pstricks}
\usepackage{wrapfig, subfigure,graphicx}
\usepackage{hyperref}
\usepackage{pdfsync}
\newcommand{\R}{{\mathbb R}} 

\newcommand{\ds}{\displaystyle}

\def \be {\mathbf{e}}

\def \p {\varphi}

\def \m{\mathbf{m}}

\renewcommand{\(}{\left(}
\renewcommand{\)}{\right)}

\newtheorem{theorem}{Theorem}[section]

\newtheorem{lemma}[theorem]{Lemma}
\newtheorem{proposition}[theorem]{Proposition}

\newtheorem{definition}[theorem]{Definition}
\newtheorem{remark}[theorem]{Remark}

\providecommand{\bysame}{\makebox[3em]{\hrulefill}\thinspace}

\begin{document}

\setcounter{equation}{0}










\title[Double metasurfaces and Optimal transport]{Double metasurfaces and Optimal transport}
\author[I. Altiner and C. E. Guti\'errez]{Irem Altiner and Cristian E. Guti\'errez}
\thanks{This work was presented at the 2024 AMS Spring Eastern Sectional Meeting at Howard University. This research was partially supported by NSF Grants DMS-1600578 and DMS-2203555.
\\ \today}
\dedicatory{Dedicated with deep admiration and warm affection to Neil Trudinger, on the occasion of his 85th birthday, in recognition of his profound and lasting contributions to mathematics.}
\address{School of Mathematical and Statistical Sciences\\Arizona State University\\Tempe, AZ 85281}
\address{Department of Mathematics\\Temple University\\Philadelphia, PA 19122}

\begin{abstract}
This paper constructs metalenses that separate homogeneous media with different refractive indices, refracting one domain into another while conserving a prescribed energy distribution. Using optimal transport theory, we design singlet and doublet metalenses satisfying energy-conserving by refraction, and obtain partial regularity of the optimal maps involved.

\end{abstract}

\maketitle

\tableofcontents

\setcounter{equation}{0}

\section{Introduction}

The purpose of this paper is to apply the theory of optimal transport (OT) to solve optical problems involving metasurfaces. Metasurfaces, or metalenses, are ultrathin optical devices built from nanostructures to manipulate light for imaging applications. They introduce
abrupt phase changes on the scale of the wavelength along the
optical path, bending light in ways not available to conventional lenses.  
By contrast, in a conventional lens the faces are chosen so that a gradual phase change accumulates as the wave
propagates through the material, reshaping the scattered wave as desired.
Metalenses are an active area of
engineering research and hold great potential for imaging applications.
They can be thinner than a sheet of paper and much lighter than
glass, and they have the potential to transform optical imaging devices---from microscopes to virtual-reality
displays and cameras, including those in smartphones; see, for example, 
\cite{science-runner-ups-2016}, \cite{2107planaroptics:capasso}, and \cite{Chen_broadband}.

Mathematically, a metalens can be modeled as a pair $(\Gamma, \Phi)$, where $\Gamma$ is a surface in three-dimensional space represented as the graph of a $C^2$ function $u$, and $\Phi$ is a $C^1$ function defined in a small neighborhood of $\Gamma$ and known as the phase discontinuity.

Optimal transport techniques have been used extensively to solve problems in free-form optics. This line of work began with reflectors in \cite{Wang:antennamasstransport}, continued with refractors in \cite{gutierrez-huang:farfieldrefractor}, and was extended to double free-form lenses in \cite{gutierrez-sabra:freeformgeneralfields}; see also \cite{2023-gutierrez:bookonOTandoptics}. 
For the design of free-form lenses, see \cite{2007-rubinstein-wolansky:freeformlens}.
Applications of OT to metasurfaces began in \cite{2020-gutierrez-merigot-thibert:metasurfacesandmasstransport}, with further results in \cite{gutierrez-sabra:chromaticaberrationinmetasurfaces}.

In this paper, we design metalenses that refract radiation to achieve prescribed energy distributions using OT techniques. 
One central problem is as follows. Suppose that light rays emanate, with prescribed directions, from a planar connected domain $\Omega_0$ located on the plane $z=0$ and situated below a given surface $\Gamma$, with intensity $\rho_0(x)$ for each $x \in \Omega_0$. 
We are also given a planar domain $\Omega_1$ on the plane $z = \beta$, above the surface $\Gamma$, together with an energy distribution specified by a density function $\rho_1$ satisfying the balance condition $\int_{\Omega_0} \rho_0(x)\, dx = \int_{\Omega_1} \rho_1(x)\, dx$.
We seek a phase discontinuity function $\Phi$ defined on $\Gamma$ such that, if $T: \Omega_0 \to \Omega_1$ denotes the map induced by the metalens $(\Gamma, \Phi)$, that is, if for $x \in \Omega_0$, $Tx \in \Omega_1$ is the image of $x$ refracted by the metalens according to the generalized Snell's law \eqref{eq:GSL}, then
\[
\int_{T^{-1}(E)}\rho_0(x)\,dx=\int_E\rho_1(x)\,dx
\]
for every measurable subset $E \subset \Omega_1$; see Figure~\ref{fig:onemetasurface}. 

Next, in Section~\ref{doublemetasurface}, we consider an analogous question involving two metalenses.
Such devices, often called \emph{doublets}, have appeared in optical engineering, notably in applications aimed at controlling chromatic aberration; see  
\cite{Arbabi:2016aa}, \cite{2022-doublet-design}, and \cite{Arbabi:2017aa}.
Specifically, let $\Omega_0$ be a planar domain on $z=0$ beneath a first metalens $(\Gamma_1, \Phi)$, and let $\Omega_1$ be a planar domain on $z= \beta$ above a second metalens $(\Gamma_2, \Psi)$. We are given density functions $\rho_i$ on $\Omega_i$, $i=0,1$, satisfying the balance condition $\int_{\Omega_0} \rho_0(x)\, dx = \int_{\Omega_1} \rho_1(x)\, dx$. We assume that the distance between the metasurfaces $\Gamma_1$ and $\Gamma_2$ is strictly positive. We seek phases $\Phi$ and $\Psi$, defined on $\Gamma_1$ and $\Gamma_2$, respectively, such that the map $T: \Omega_0 \to \Omega_1$ induced by these phases in accordance with the generalized Snell's law \eqref{eq:GSL} satisfies
\[
\int_{T^{-1}(E)}\rho_0(x)\,dx=\int_E\rho_1(x)\,dx
\]
for every measurable subset $E \subset \Omega_1$; see Figure~\ref{fig:two metasurfaces}. 

We solve these problems by introducing cost functions naturally associated with the optical setups, formulas \eqref{eq:cost for meta on z=f(x,y)3} and \eqref{eq:cost for double meta}, and applying optimal transport theory.   
The construction of the phases follows from the corresponding optimal maps, under suitable assumptions on the cost functions.
Moreover, under additional assumptions on the surfaces entering the costs, we prove partial regularity of the optimal maps; this analysis is carried out in Section~\ref{sec:regularity of the optimal map and costs}.   

The paper is organized as follows. 
Section~\ref{sec:preliminaries} introduces the generalized Snell laws of refraction and reflection.
Section~\ref{sec:optimal transport} presents the relevant results from optimal transport theory.
Section~\ref{sec:construction of metasurfaces} contains the construction of phases for a single metasurface (Section~\ref{sec:single metasurface}, Theorem~\ref{thm:result for one metasurface}) and for two metasurfaces (Section~\ref{doublemetasurface}, Theorem~\ref{thm:main theorem for two metasurfaces}).
In Section~\ref{sec:regularity of the optimal map and costs}, we analyze the geometric and analytic conditions on the surfaces $\Gamma_1,\Gamma_2$ that allow the application of the partial regularity results for optimal maps from \cite{partial-regularity_optimalmaps} to the cost~\eqref{eq:general formula for the cost}; see Theorem~\ref{thm:main result of regularity for general costs}.

We would like to thank Alessio Figalli for his comments and suggestions on an earlier version of this paper \cite{2025:arxivfirstversion}. His feedback helped clarify the presentation of the results in Section~\ref{sec:construction of metasurfaces}, which were established in that earlier version using flows of vector fields in combination with the partial regularity theory discussed in Section~\ref{sec:regularity of the optimal map and costs}.

\setcounter{equation}{0}

\section{Preliminaries}\label{sec:preliminaries}

\subsection{Generalized Snell's law}


We begin by recalling the generalized Snell's law.

Let $\Gamma$ be a surface in $\mathbb{R}^3$ that separates two media, $I$ and $II$, with refractive indices $n_1$ and $n_2$, respectively. Suppose that $\Gamma$ is defined by the equation $\phi(x,y,z) = 0$, where $\phi \in C^1$. Let $\psi$ be a $C^1$ function defined in a small neighborhood of $\Gamma$. Given points $A \in I$ and $B \in II$, we consider paths from $A$ to $B$ that pass through a point $P$ on $\Gamma$, and we seek the point $P$ that minimizes the total travel time.

The velocities of propagation in $I$ and $II$ are $v_1 = c/n_1$ and $v_2 = c/n_2$, respectively, where $c$ is the speed of light in vacuum. If $-\psi(P)$ represents the height of the obstacle at point $P \in \Gamma$, then the total travel time from $A$ to $B$, passing through $P$, is given by
\[
\frac{n_1}{c}|A-P|+\frac{n_2}{c}|B-P|-\frac{1}{c}\psi(P).
\]
The objective is to minimize this expression over $P \in \Gamma$.
Multiplying by $c$ and applying Fermat's principle, we find that $P$ is a critical point of
$$n_1\, |P-A|+n_2\,|B-P| -\psi(P),$$
subject to the constraint $\phi(P) = 0$ for $P \in \Gamma$.
Using the method of Lagrange multipliers, we obtain 
$$\nabla\left(n_1\, |P-A|+n_2\,|B-P| -\psi(P)\right)\times \nabla \phi(P)=0.$$
Since $\nabla \phi$ is parallel to the normal vector $\nu$ to $\Gamma$, and
$$\nabla\left(n_1\, |P-A|+n_2\,|B-P| -\psi(P)\right)=n_1\dfrac{{P-A}}{|P-A|}-n_2\dfrac{{B-P}}{|B-P|}-\nabla \psi(P),$$ 
denoting the unit directions of the incident and refracted (or transmitted) rays by ${\bf x}=\dfrac{{P-A}}{|P-A|}$ and ${\bf m}=\dfrac{{B-P}}{|B-P|}$, respectively, we obtain the generalized Snell's law of refraction (GSL)
\begin{equation*}
\left(n_1\,{\bf x}-n_2\,{\bf m}\right)\times \nu=\nabla \psi\times \nu.
\end{equation*}
This is equivalent to
\begin{equation}\label{eq:GSL}
n_1\,{\bf x}-n_2\,{\bf m}=\lambda \,\nu+\nabla \psi,
 \end{equation}
where $\lambda\in \R$ can be explicitly calculated in terms of $n_1, n_2, {\bf x},\nu$, and $\nabla \psi$; see \cite{gps} and \cite{gutierrez-sabra:chromaticaberrationinmetasurfaces}.
 
We can also deduce the generalized Snell's law of reflection. Indeed, taking both points $A$ and $B$ in medium $I$ and proceeding in a similar way, we obtain 
\begin{equation}\label{eq:GSLreflection}
n_1{\bf x} -  \,n_1{\bf r}=\lambda \,\nu+\nabla \psi,
\end{equation}
for some $\lambda\in \R$, where now ${\bf r}$ represents the unit direction of the reflected ray.

\subsection{Optimal Transport}\label{sec:optimal transport}
In this section, we briefly review a few facts from optimal transport theory that will be used later. 
We follow the approach from \cite[Chapter 6]{2023-gutierrez:bookonOTandoptics}.

Let $(D,d)$ and $(D^*,d^*)$ be compact metric spaces, and let $c:D\times D^*\to \R$ be a Lipschitz cost function; that is, there exists a constant $K>0$ such that 
\[
|c(x_1,y_1)-c(x_2,y_2)|\leq K\,\(d(x_1,x_2)+d^*(y_1,y_2)\),
\] 
for all $x_i\in D$, $y_i\in D^*$, $i=1,2$.

Given $u\in C(D)$, the $c$-transform of $u$ is by definition
\[
u^c(m)=\inf_{x\in D}\left\{c(x,m)-u(x)\right\}
\]
for each $m\in D^*$.

\begin{definition}
The function $\phi:D\to \R$ is $c$-concave if, for each $x_0\in D$, there exist $m_0\in D^*$ and $b\in \R$ such that $\phi(x)\leq c(x,m_0)-b$ for all $x\in D$, with equality at $x=x_0$.
\end{definition}
Notice that each $c$-concave function is Lipschitz in $D$.

\begin{definition}
Given $\phi:D\to\R$, the $c$-superdifferential of $\phi$ 
is the multivalued map from $D$ to $\mathcal P(D^*)$ defined by
\[
\partial_c\phi(x)=\left\{m\in D^*: \phi(x)+\phi^c(m)=c(x,m) \right\}.
\] 
\end{definition} 
Notice that if $\phi$ is $c$-concave, then $\partial_c\phi(x)\neq \emptyset$ for all $x\in D$.

Let $\mu$ be a Borel measure on $D$.
We introduce the following assumption on the cost $c$:

\indent \paragraph{\bf SV}\label{single valued a.e.}
\textit{For each $c$-concave function $\phi:D\to \R$, the set
\[
\left\{x\in D: \partial_c\phi(x) \text{ is not a singleton}\right\}
\]
has $\mu$-measure zero.\footnote{For an analogue of this condition in the multi-marginal Monge problem, with an application to metasurfaces, see \cite{altiner-gutierrez:multimarginalmongeproblemmetasurfaces}.}}

We have the following lemma.
\begin{lemma}
If $c$ satisfies \nameref{single valued a.e.}, then $\partial_c \phi$ is measurable for each $c$-concave function $\phi$; that is, for each Borel subset $F\subset D^*$, the set 
$(\partial_c \phi)^{-1}(F)=\left\{ x\in D: \partial_c\phi(x)\cap F\neq \emptyset \right\}$ is a $\mu$-measurable subset of $D$.
\end{lemma}
\begin{proof}
We first show that 
\[
\mathcal C=\left\{E\subset D^*: (\partial_c \phi)^{-1}(E)\text{ is a $\mu$-measurable subset of $D$} \right\}
\]
is a $\sigma$-algebra. Since $\phi$ is $c$-concave, $\partial_c\phi(x)\neq \emptyset$ for each $x\in D$, and so $(\partial_c \phi)^{-1}(D^*)=D$, i.e., $D^*\in \mathcal C$. The class $\mathcal C$ is clearly closed under countable unions. It is also closed under complements, since
\[
(\partial_c \phi)^{-1}(E^c)=\((\partial_c \phi)^{-1}(E)\)^c\cup \((\partial_c \phi)^{-1}(E^c)\cap (\partial_c \phi)^{-1}(E)\)
\] 
where the last term in the union has measure zero by \nameref{single valued a.e.}.
We next show that $\mathcal C$ contains all compact subsets of $D^*$.
Let $K\subset D^*$ be compact and let $\{u_j\}_{j=1}^\infty$ be a sequence in $(\partial_c \phi)^{-1}(K)\subset D$. 
Since $D$ is compact, after passing to a subsequence we may assume that $u_j\to u^0$ for some $u^0\in D$.
Since $u_j\in (\partial_c \phi)^{-1}(K)$, there exists $m_j\in \partial_c\phi(u_j)\cap K$ for each $j$. Since $K$ is compact, after passing to a further subsequence we have $m_{j_\ell}\to m_0\in K$.
We then have 
$
\phi(u_{j_\ell})+\phi^c(m_{j_\ell})=
c\(u_{j_\ell},m_{j_\ell}\)
$
and letting $\ell\to \infty$, using the continuity of $\phi$, $\phi^c$, and $c$, yields 
$
\phi(u^0)+\phi^c(m_0)=
c\(u^0,m_0\)
$
that is, $m_0\in \partial_c\phi(u^0)$ and hence $u^0\in (\partial_c\phi)^{-1}(K)$.
Therefore $(\partial_c\phi)^{-1}(K)$ is closed in $D$, and so it is compact and $\mu$-measurable. Thus $K\in\mathcal C$ for every compact $K\subset D^*$. Since $D^*$ is a compact metric space, its Borel sets are generated by compact sets; hence $\mathcal C$ contains every Borel subset of $D^*$.
\end{proof}


\begin{definition}
Let $s:D\to \mathcal P(D^*)$ be a multivalued map such that the set 
\[
\left\{x\in D: s(x) \text{ is not a singleton}\right\}
\] 
has $\mu$-measure zero and $s$ is a measurable map.
If $\mu^*$ is a Borel measure on $D^*$, we say that $s$ is measure preserving from $(D,\mu)$ to $(D^*,\mu^*)$ if 
\[
s_\# \mu(F):=\mu(s^{-1}(F))=\mu^*(F)
\]
for each Borel subset $F\subset D^*$. We denote by $\mathcal S(\mu,\mu^*)$ the class of all such measure-preserving maps.
\end{definition}
Then $s_\#\mu$ is a Borel measure on $D^*$; see \cite[Chapter 5]{2023-gutierrez:bookonOTandoptics}.

If $\mu(D)=\mu^*(D^*)$, then recall that the Monge problem is to minimize 
\[
\int_D c\(x,s(x)\)\,d\mu(x),
\]
among all maps $s\in \mathcal S(\mu,\mu^*)$.
This problem has a solution and it is unique under the following circumstances.

\begin{theorem}\label{thm:existence and uniqueness of optimal maps}
If $\mu(D)=\mu^*(D^*)$ and the cost $c$ satisfies condition \nameref{single valued a.e.}, then there exists a $c$-concave function $\phi$ such that the minimum of the Monge problem is attained at $s=\partial_c\phi$.
In addition, if $\mu(G)>0$ for every nonempty open set $G\subset D$, then the minimizer is unique.

\end{theorem}
For a proof of this theorem, see \cite[Lemmas 6.6 and 6.7]{2023-gutierrez:bookonOTandoptics}.

\subsubsection{Application to $\R^n$}
We first prove the following lemma.
\begin{lemma}\label{lm:injectivity implies A and B}
Suppose $D=\Omega_0$ and $D^*=\Omega_1$ are compact domains in $\R^n$, $|\partial \Omega_0|=0$, and the cost $c:\Omega_0\times \Omega_1\to\R$ is $C^1$.
If, for each $x\in \Omega_0$, the map from $\Omega_1$ to $\R^n$ given by $y\mapsto \nabla_xc(x,y)$ is injective, then $c$ satisfies condition \nameref{single valued a.e.} with $\mu$ equal to Lebesgue measure.
\end{lemma}
\begin{proof}
Let $\phi:\Omega_0\to \R$ be $c$-concave. Then $\phi$ is Lipschitz in $\Omega_0$, and by Rademacher's theorem there is a set $N\subset \Omega_0$ of Lebesgue measure zero such that $\phi$ is differentiable in $\Omega_0\setminus N$.
Let $x_1\in \Omega_0\setminus \(N\cup \partial \Omega_0\)$ and let $m_1,m_2\in \partial_c\phi (x_1)$.
Then $\phi(x_1)+\phi^c(m_i)=c(x_1,m_i)$ for $i=1,2$.
This means $\phi(x_1)+\inf_{z\in \Omega_0}\(c(z,m_i)-\phi(z)\)=c(x_1,m_i)$, that is,
\[
\phi(x_1)+c(z,m_i)-\phi(z)\geq c(x_1,m_i)\quad \forall z\in \Omega_0
\]
or equivalently 
\[
\phi(x_1)-c(x_1,m_i)\geq \phi(z)-c(z,m_i)\quad \forall z\in \Omega_0
\]
and so the maximum of $\phi(z)-c(z,m_i)$ over $\Omega_0$ is attained at $z=x_1$.
Since $c$ is $C^1$ and $\phi$ is differentiable at $x_1$, we get that $\nabla_z \phi(z)-\nabla_z c(z,m_i)=0$ at $z=x_1$ for $i=1,2$. The injectivity of $\nabla_x c(x_1,\cdot)$ then implies $m_1=m_2$.
In particular, at such points we have $\nabla \phi(x_1)-\nabla_x c\(x_1,\partial_c\phi (x_1)\)=0$.
\end{proof}

From Theorem~\ref{thm:existence and uniqueness of optimal maps}, we obtain the following proposition, which will be used later.

\begin{proposition}\label{prop:representation formula from injectivity}
Suppose $\Omega_0,\Omega_1$ are compact domains in $\R^n$, $\mu$ is a Borel measure on $\Omega_0$ satisfying $\mu(E)=\int_E g(x)\,dx$ for every Borel set $E\subset \Omega_0$ and some $g>0$ a.e. in $\Omega_0$, $\mu^*$ is a Borel measure on $\Omega_1$, $\mu(\Omega_0)=\mu^*(\Omega_1)$, and $|\partial \Omega_0|=0$.

If the cost function $c$ is $C^1$ and the map $y\mapsto \nabla_x c(x,y)$ is injective for each $x\in \Omega_0$,
then there exists a $c$-concave function $\psi$ such that the unique optimal map $T=\partial_c\psi$ for the Monge problem satisfies
\begin{equation}\label{eq:representation formula for gradient of c}
\nabla \psi(x)=\nabla_x c\(x,Tx\)
\end{equation}
for a.e. $x\in \Omega_0$.
\end{proposition}
\begin{proof}
From the form of $\mu$ we have $\mu(N)=0$ if and only if $|N|=0$. Hence, by Lemma~\ref{lm:injectivity implies A and B},
the cost $c$ satisfies \nameref{single valued a.e.} with the measure $\mu$. 
Moreover, since $g>0$ a.e. in $\Omega_0$, we have $\mu(G)>0$ for every nonempty open set $G\subset \Omega_0$. Therefore, by Theorem~\ref{thm:existence and uniqueness of optimal maps}, the optimal map is unique and is given by $T=\partial_c\psi$ for some $c$-concave function $\psi$.
Formula \eqref{eq:representation formula for gradient of c} follows from the proof of Lemma~
\ref{lm:injectivity implies A and B}.

\end{proof}

\begin{remark}
For extensions of the results of this section to the multi-marginal Monge problem, see \cite{altiner-gutierrez:multimarginalmongeproblemmetasurfaces}.
\end{remark}

\setcounter{equation}{0}
\section{Construction of metasurfaces}\label{sec:construction of metasurfaces}

In this section, we use the results from Section~\ref{sec:optimal transport} to construct metasurfaces that refract radiation in a prescribed manner in several geometric scenarios. We begin with the simpler case of a single metasurface.

\subsection{Single metasurfaces}\label{sec:single metasurface}
A planar domain $\Omega_0$ is given, and for each $x=(x_1,x_2) \in \Omega_0$, a unit field of directions 
$\be(x) = \(e_1(x), e_2(x), e_3(x)\)$ is prescribed, satisfying $e_3(x) > 0$. A surface $S_1$ above the plane $z = 0$ is given as the graph of a function $f$ and is related to the field $\be$ as follows. 
For each $x \in \Omega_0$, the ray $\{(x,0)+t\be(x): t>0\}$ intersects the graph of $f$ at a unique point $\(\p(x), f\(\p(x)\)\) \in S_1$. Here, $\p: \Omega_0 \to \Omega_0'$ is a one-to-one $C^2$ map onto another planar domain $\Omega_0'$, and $f \in C^1\(\Omega_0'\)$.
In the collimated case, that is, when $\be(x)=(0,0,1)$, it is clear that $\Omega_0'=\Omega_0$ and $\p$ is the identity.

We are also given a plane $z = \beta$ located at a positive distance above $S_1$, a second domain $\Omega_1 \subset \R^2$, and densities $\rho_0$ on $\Omega_0$ and $\rho_1$ on $\Omega_1$ satisfying the global energy conservation condition:
\begin{equation}\label{eq:conservation of energy global}
\int_{\Omega_0} \rho_0(x)\, dx = \int_{\Omega_1} \rho_1(x)\, dx.
\end{equation}
A material with refractive index $n_1$ fills the region between the plane $z = 0$ and the surface $S_1$, while a material with refractive index $n_2$ fills the region between the surface $S_1$ and the plane $z = \beta$.

We consider the following problem for metalenses. A light ray emitted from $x \in \Omega_0$ in the direction $\be(x)$ strikes $S_1$ at the point $\(\p(x), f(\p(x))\)$. 
On the surface $S_1$, a phase discontinuity function $\Phi$ is defined so that the ray is refracted into a unit direction $\m(x) = (m_1(x), m_2(x), m_3(x))$ with $m_3(x) > 0$, in accordance with the generalized Snell law \eqref{eq:GSL}, which depends on the phase $\Phi$, and proceeds to a point $(Tx, \beta)$ on the plane $z = \beta$.
That is, each point $x \in \Omega_0$ is mapped to a point $(Tx, \beta)$.

\begin{figure}[t]
    \centering
    \includegraphics[scale=0.45]{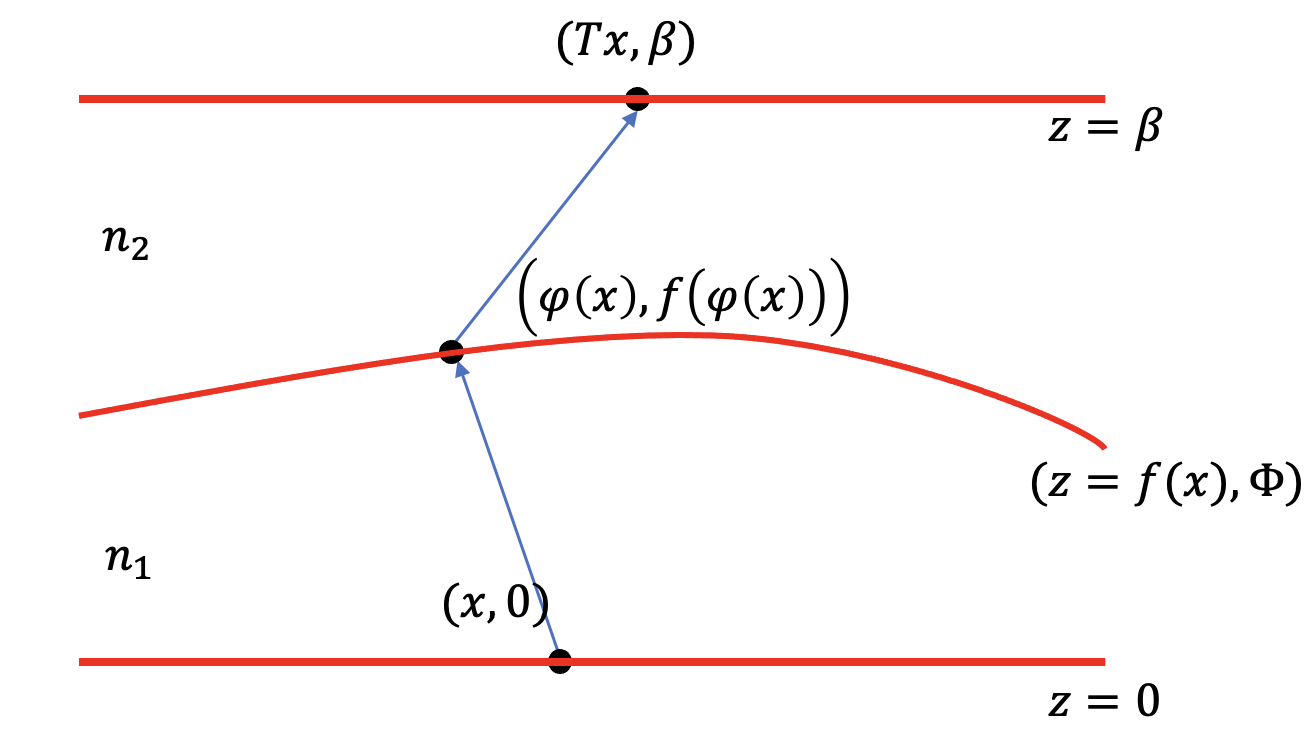}
    \caption{Single metasurface}
    \label{fig:onemetasurface}
\end{figure}

The goal is to determine a phase function $\Phi$ on the surface $z = f(x)$ such that it is tangential to the surface, i.e.,
\begin{equation}\label{eq:tangential condition}
\nabla \Phi\(x, f(x)\) \cdot (-\nabla f(x), 1) = 0 \quad \text{for all } x \in \Omega_0',
\end{equation}
and such that the mapping $T$ satisfies $T\(\Omega_0\) = \Omega_1$ and the local conservation of energy condition
\begin{equation}\label{eq:conservation of energy}
\int_{T^{-1}(E)} \rho_0(x)\, dx = \int_E \rho_1(x)\, dx
\end{equation}
for every Borel set $E \subset \Omega_1$; see Figure~\ref{fig:onemetasurface}.

We now explain the tangential condition \eqref{eq:tangential condition}. The phase is initially defined only on
$
S_1=\{X(x):=(x,f(x)):x\in\Omega_0'\}\subset\mathbb R^3.
$
To compute its ambient gradient, one must extend it to a tubular neighborhood of $S_1$. There are many possible extensions, but the tangentiality requirement selects the extension that is constant along normal lines. Let $\Phi:S_1\to\mathbb R$ be a smooth phase function, and write $\phi(x):=\Phi(X(x))$ for $x\in\Omega_0'.$ The unit normal to $S_1$ is
\[
\nu(x)=\frac{N(x)}{\sqrt{1+|\nabla f(x)|^2}},
\qquad N(x)=(-\nabla f(x),1).
\]
If $U_1$ is a sufficiently small neighborhood of $S_1$, then every $Y\in U_1$ can be written uniquely as
\[
Y=X(\xi)+t\nu(\xi),
\qquad \xi\in\Omega_0',
\]
with $|t|$ small. Using the same notation for the extension, define
\[
\Phi\bigl(X(\xi)+t\nu(\xi)\bigr)=\phi(\xi).
\]
Thus $\Phi$ is constant along each normal line. Differentiating this identity with respect to $t$ at $t=0$ gives
\[
\nabla\Phi(X(x))\cdot\nu(x)=0,
\]
which is equivalent to \eqref{eq:tangential condition}. With this convention, $\nabla\Phi$ along $S_1$ is the intrinsic surface gradient of the phase.

To find the phase $\Phi$, we use optimal transport theory with an appropriate cost depending on the surface $S_1$ and the plane $z=\beta$.

Let us first analyze the trajectory of the ray.
The light ray starts from the point $(x,0)$ with direction $\be(x)$ and travels to $\(\p(x),f(\p(x))\)\in S_1$; from this point it travels to $(Tx,\beta)$. Thus the refracted unit direction is 
\begin{equation}\label{eq:formula for m general field non planar}
\m(x)=\dfrac{\(Tx-\p(x),\beta-f(\p(x))\)}{\sqrt{|Tx-\p(x)|^2+\(\beta-f(\p(x))\)^2}}.
\end{equation}

On the other hand, according to the GSL \eqref{eq:GSL}, the phase $\Phi$ must satisfy on $S_1$
\begin{equation}\label{eq:GSL for one metasurface}
    n_1\textbf{e}(x) - n_2\textbf{m}(x) = \lambda\nu(\p(x),f(\p(x))) + \nabla\Phi(\p(x),f(\p(x))) 
\end{equation}
for some $\lambda\in \R$, where $\nu(\p(x),f(\p(x))) = (-\nabla f(\p(x)),1)$ is the normal vector to the surface $z=f(x)$ at the point $\(\p(x),f(\p(x))\)$. Since $\Phi$ is required to be tangential, we have
\begin{equation}\label{eq:tangential condition single}
-\Phi_{x_1}(\p(x),f(\p(x)))\,f_{x_1}(\p(x))-\Phi_{x_2}(\p(x),f(\p(x)))\,f_{x_2}(\p(x))+\Phi_{x_3}(\p(x),f(\p(x)))=0.
\end{equation}
From \eqref{eq:GSL for one metasurface},
\begin{equation}\label{eq:general e m nonplanar single}
\begin{cases}
n_1\,e_1(x)-n_2\,m_1(x)&=-\lambda\,f_{x_1}(\p(x))+\Phi_{x_1}(\p(x),f(\p(x)))\\
n_1\,e_2(x)-n_2\,m_2(x)&=-\lambda\,f_{x_2}(\p(x))+\Phi_{x_2}(\p(x),f(\p(x)))\\
n_1\,e_3(x)-n_2\,m_3(x)&=\lambda+\Phi_{x_3}(\p(x),f(\p(x))).
\end{cases}
\end{equation}
%
Thus from \eqref{eq:formula for m general field non planar}\begin{align}
\label{eq:formula for Tx-phix}
Tx-\p(x)&=\sqrt{(\beta-f(\p(x)))^2+|Tx-\p(x)|^2}\\
&\qquad \(\frac{n_1}{n_2}\(e_1(x),e_2(x)\)+\dfrac{1}{n_2}\,
\(\lambda\,\nabla f(\p(x))-\(\Phi_{x_1}\(\p(x),f(\p(x))\),\Phi_{x_2}\(\p(x),f(\p(x))\)\)\)\).\notag
\end{align}
The last equation in \eqref{eq:general e m nonplanar single} and \eqref{eq:tangential condition single} yield 
\begin{align*}
\lambda&=n_1\,e_3(x)-n_2\,\dfrac{\beta-f(\p(x))}{\sqrt{(\beta-f(\p(x)))^2+|Tx-\p(x)|^2}}-\Phi_{x_3}\\
&=n_1\,e_3(x)-n_2\,\dfrac{\beta-f(\p(x))}{\sqrt{(\beta-f(\p(x)))^2+|Tx-\p(x)|^2}}-\Phi_{x_1}\,f_{x_1}-\Phi_{x_2}\,f_{x_2},
\end{align*}
where $\Phi_{x_i}$ is evaluated at $(\p(x), f(\p(x)))$, and $f_{x_i}$ is evaluated at $\p(x)$, $x\in\Omega_0$. To simplify the notation, we omit these points in the calculations below.
Hence,
\[
\lambda\,f_{x_i}-\Phi_{x_i}=
\(n_1e_3(x)-n_2\,\dfrac{\beta-f(\p(x))}{\sqrt{(\beta-f(\p(x)))^2+|Tx-\p(x)|^2}}-\Phi_{x_1}\,f_{x_1}-\Phi_{x_2}\,f_{x_2}\)\,f_{x_i}-\Phi_{x_i},\]
for $i=1,2$, and so 
\begin{align*}
&\lambda\,\nabla f-\(\Phi_{x_1},\Phi_{x_2}\)\\
&=
\(n_1e_3(x)-n_2\,\dfrac{\beta-f(\p(x))}{\sqrt{(\beta-f(\p(x)))^2+|Tx-\p(x)|^2}}\)\nabla f
-
\(\nabla f\otimes \nabla f\)\(\Phi_{x_1},\Phi_{x_2}\)-\(\Phi_{x_1},\Phi_{x_2}\)\\
&=
n_1e_3(x)\,\nabla f
-n_2\,\dfrac{\beta-f(\p(x))}{\sqrt{(\beta-f(\p(x)))^2+|Tx-\p(x)|^2}}\,\nabla f
-
\(Id+\nabla f\otimes \nabla f\)\(\Phi_{x_1},\Phi_{x_2}\).
\end{align*}
Substituting the last expression into \eqref{eq:formula for Tx-phix} shows that $T$ and $\Phi$ are related by the equation
\begin{equation}\label{eq:main formula for T2 non flat}
\begin{split}
\dfrac{Tx-\p(x)+(\beta-f(\p(x)))\nabla f(\p(x))}{\sqrt{(\beta-f(\p(x)))^2+|Tx-\p(x)|^2}}
&=\dfrac{-1}{n_2}\(Id+\nabla f\otimes \nabla f\)\(\Phi_{x_1},\Phi_{x_2}\)\\
&
+\dfrac{n_1}{n_2}e_3\nabla f+\dfrac{n_1}{n_2}(e_1(x),e_2(x)).
\end{split}
\end{equation}
We now connect this with the optimal transport theory.
Let us introduce the cost function
\begin{equation}\label{eq:cost for meta on z=f(x,y)3}
c(x,y)=\sqrt{(\beta-f(\p(x)))^2+|\p(x)-y|^2},
\end{equation}
where $(x,y) \in \Omega_0 \times \Omega_1$.

Notice that if $J_\p=\left(\dfrac{\partial \p_i}{\partial x_j}\right)_{i,j=1}^2$ 
denotes the Jacobian matrix of $\p$, then from \eqref{eq:cost for meta on z=f(x,y)3}
we get that 
\begin{equation}\label{eq:formula for Dx of cost2}
\nabla_xc(x,y)=J_\p(x)\,\dfrac{\p(x)-y-(\beta-f(\p(x)))\nabla f(\p(x))}{\sqrt{(\beta-f(\p(x)))^2+|\p(x)-y|^2}}.
\end{equation}


Since the map $\p$ is invertible, using \eqref{eq:formula for Dx of cost2} we can rewrite \eqref{eq:main formula for T2 non flat}
as follows:
\begin{equation}\label{eq:formula for Dx of cost optimal map3}
\(J_\p(x)\)^{-1}\nabla_xc(x,Tx)=\dfrac{1}{n_2}\(Id+\nabla f\otimes \nabla f\)\(\Phi_{x_1},\Phi_{x_2}\)
-\dfrac{n_1}{n_2}e_3(x)\nabla f - \dfrac{n_1}{n_2}(e_1(x),e_2(x)),
\end{equation}
showing that $T$ and $\Phi$ are related via the cost $c$.
At this point, we do not yet know whether the map $T$ satisfies the conservation condition \eqref{eq:conservation of energy}.

If $T$ is the optimal transport map with respect to the cost $c$ and densities $\rho_0,\rho_1$,
then $T$ satisfies the energy conservation condition \eqref{eq:conservation of energy}, and the phase $\Phi$ is given by \eqref{eq:formula for Dx of cost optimal map3}. 
Moreover, if the assumptions of Proposition~\ref{prop:representation formula from injectivity} hold, then 
by \eqref{eq:representation formula for gradient of c}, $\nabla_xc(x,Tx)$ equals the gradient of a $c$-concave function $\psi$ a.e., and therefore \eqref{eq:formula for Dx of cost optimal map3} becomes
\begin{equation}\label{eq:formula for Dx of cost optimal map3BIS}
\(J_\p(x)\)^{-1}\nabla\psi(x)=\dfrac{1}{n_2}\(Id+\nabla f\otimes \nabla f\)\(\Phi_{x_1},\Phi_{x_2}\)
-\dfrac{n_1}{n_2}e_3(x)\nabla f - \dfrac{n_1}{n_2}(e_1(x),e_2(x)),
\end{equation}
and so the phase $\Phi$ is determined by this equation.
 
Moreover, we can solve equation \eqref{eq:formula for Dx of cost optimal map3BIS}
for $\(\Phi_{x_1},\Phi_{x_2}\)$. Indeed,
from \eqref{eq:formula for Dx of cost optimal map3BIS}
\[
n_2\,\(J_\p(x)\)^{-1}\nabla\psi(x)+n_1 e_3(x) \nabla f + n_1\,(e_1(x),e_2(x))
=\(Id+\nabla f\otimes \nabla f\)\(\Phi_{x_1},\Phi_{x_2}\),
\] 
Once $T,c, \p$, and $f$ are known, it remains to solve this equation for $\(\Phi_{x_1},\Phi_{x_2}\)$.
Recall the Sherman--Morrison formula: if $A$ is an $n\times n$ invertible matrix, $u,v$ are $n$-column vectors, $u\otimes v=uv^t$, and $1+v^tA^{-1}u\neq 0$, then $A+u\otimes v$ is invertible and
\begin{equation}\label{eq:Sherman Morrinson formula}
\(A+u\otimes v\)^{-1}=A^{-1}-\dfrac{A^{-1}u\otimes vA^{-1}}{1+v^tA^{-1}u},\qquad \det\(A+u\otimes  v\)=(1+u^t A^{-1} v)\,\det A.
\end{equation}
Applying this formula with $A=Id$ and $u=v=\nabla f$, we get 
\[
\(Id+\nabla f\otimes \nabla f\)^{-1}
=
Id-\dfrac{\nabla f\otimes \nabla f}{1+|\nabla f|^2}.
\]
Therefore, the phase $\Phi$ is given by  
\begin{equation}\label{eq:formula for the first phase Phi}
\(\Phi_{x_1},\Phi_{x_2}\)
=
\(Id-\dfrac{\nabla f\otimes \nabla f}{1+|\nabla f|^2} \)\(n_2\,\(J_\p(x)\)^{-1}\nabla\psi(x)+n_1 e_3(x)\, \nabla f + n_1\,(e_1(x),e_2(x))\),
\end{equation}
where $\Phi_{x_i}$ is evaluated at $\(\p(x),f(\p(x))\)$, and $\nabla f$ is evaluated at $\p(x)$, $x\in \Omega_0$.

In the particularly important case when $\be(x)=(0,0,1)$ and $f$ is constant, we obtain the formula
\begin{equation}\label{eq:e is vertical and f is constant}
\(\Phi_{x_1},\Phi_{x_2}\)
=
n_2\,\nabla\psi(x).
\end{equation}

 
We summarize this in the following theorem. 

\begin{theorem}\label{thm:result for one metasurface}
With the setup from the beginning of this section, suppose that $f\in C^1(\Omega_0')$, the cost $c$ given by \eqref{eq:cost for meta on z=f(x,y)3} satisfies the assumptions of Proposition~\ref{prop:representation formula from injectivity}\footnote{The injectivity of the map $y\mapsto \nabla_x c(x,y)$ holds under assumption \eqref{eq:C1condition for beta in terms of max and min}; see Section~\ref{subsect:verification of twist g=beta}.}, and \eqref{eq:conservation of energy global} holds with  
$\rho_0$ strictly positive a.e.,
then there is a $c$-concave function $\psi$ such that the phase $\Phi$ satisfies
\[
\(\Phi_{x_1},\Phi_{x_2}\)
=
\(Id-\dfrac{\nabla f\otimes \nabla f}{1+|\nabla f|^2} \)\(n_2\,\(J_\p(x)\)^{-1}\nabla\psi(x)+n_1 e_3(x)\, \nabla f + n_1\,(e_1(x),e_2(x))\),
\]
where $\psi$ is the function in \eqref{eq:representation formula for gradient of c}, with $\Phi_{x_i}$ evaluated at $(\varphi(x), f(\varphi(x)))$, and $\nabla f$ evaluated at $\varphi(x)$.

\end{theorem}

In other words, the phase $\Phi$ can be obtained from the optimal map $T$ with respect to the cost \eqref{eq:cost for meta on z=f(x,y)3} and the densities $\rho_0,\rho_1$ through the representation formula 
\eqref{eq:formula for the first phase Phi}.
In the collimated case, i.e., when $\be=(0,0,1)$ and $f$ is constant, \eqref{eq:e is vertical and f is constant} shows that $\Phi$ is equal, up to a constant, to $n_2\psi$.

\subsection{Double metasurfaces} \label{doublemetasurface}

As in Section~\ref{sec:single metasurface}, we are given a planar domain $\Omega_0$, and for each $x = (x_1, x_2) \in \Omega_0$, a unit field of directions 
$\be(x) = \(e_1(x), e_2(x), e_3(x)\)$ is prescribed, satisfying $e_3(x) > 0$. 
A surface $S_1$ above the plane $z = 0$ is given as the graph of a strictly positive function $f$, which is related to the field $\be$ as in Section~\ref{sec:single metasurface}; that is,  
for each $x \in \Omega_0$, the ray in the direction $\be(x)$ intersects the graph of $f$ at a unique point $\(\p(x), f\(\p(x)\)\) \in S_1$, where $\p: \Omega_0 \to \Omega_0'$ is a $C^2$ one-to-one mapping onto another planar domain $\Omega_0'$, and $f \in C^1\(\Omega_0'\)$. 

Next, we are given a surface $S_2$ described by the graph of a function $g: \Omega_1 \to \R^+$ such that the distance between $S_1$ and $S_2$ is strictly positive, where $\Omega_1$ is a planar domain. 
The region between the plane $z = 0$ and the graph of $f$ is filled with a material of refractive index $n_1$; the region between the graphs of $f$ and $g$ has refractive index $n_2$; and the region above the graph of $g$ is filled with a material of refractive index $n_3$. 
We are also given a plane $z = \beta$ located above the graph of $g$. 

We now seek two phases: $\Phi$ on the surface $S_1$, tangential to $S_1$, and $\Psi$ on the surface $S_2$, tangential to $S_2$, that solve the following problem. 
For each $x = (x_1, x_2) \in \Omega_0$, a ray emitted from $x$ in the direction $\be(x)$ strikes $S_1$ at the point $\(\p(x), f(\p(x))\)$ and is refracted into a unit direction $\m(x)$ according to the generalized Snell's law (GSL) with respect to the phase $\Phi$, reaching the point $(Tx, g(Tx))$ on $S_2$. 
At $S_2$, the ray is then refracted into the vertical direction $(0, 0, 1)$ according to the GSL, now with respect to the phase $\Psi$, and continues until it reaches the plane $z = \beta$ at the point $(Tx, \beta)$ in the domain $\Omega_1 \times \{\beta\}$. 

\begin{figure}[t]
    \centering
    \includegraphics[scale=0.45]{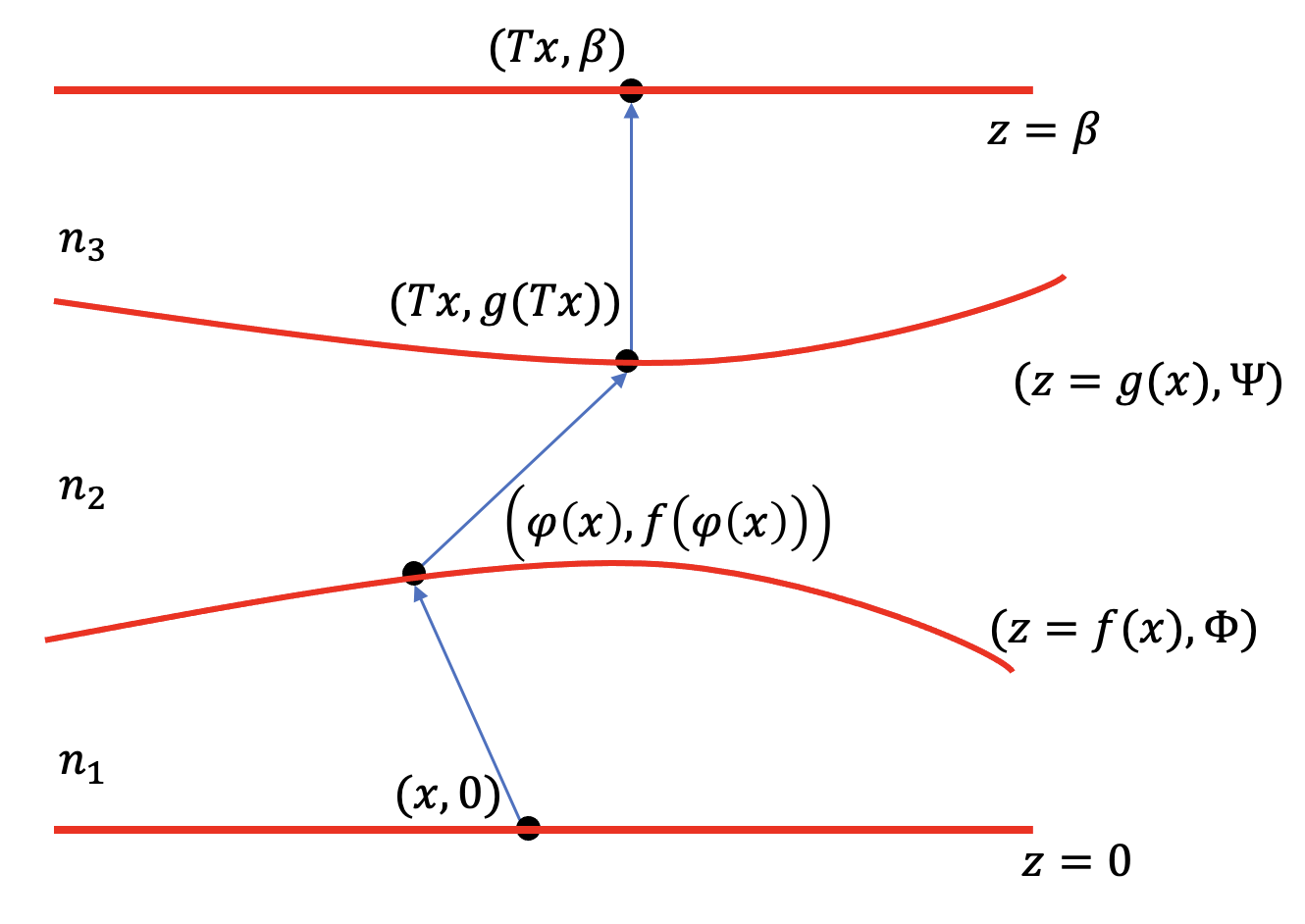}
    \caption{Double metasurface}
\label{fig:two metasurfaces}
\end{figure}

Given densities $\rho_0$ on $\Omega_0$ and $\rho_1$ on $\Omega_1$ satisfying \eqref{eq:conservation of energy global}, we seek phases $\Phi$ and $\Psi$ such that the map $T: \Omega_0 \to \Omega_1$ satisfies \eqref{eq:conservation of energy}; see Figure~\ref{fig:two metasurfaces}.

As before, to find the phases $\Phi$ and $\Psi$, we use optimal transport theory with an appropriate cost depending now on the surfaces $S_1$ and $S_2$.

Let us first analyze the trajectory of the ray.
The light ray starts from the point $(x,0)$ with unit direction $\be(x)$ and travels to $\(\p(x),f(\p(x))\)\in S_1$; from this point it travels to $(Tx,g(Tx))$ with unit direction $\m(x)$. Finally, the ray travels vertically, i.e., with direction $(0,0,1)$, to the final point $(Tx, \beta)$. Thus the refracted unit direction between $S_1$ and $S_2$ is 
\begin{equation}\label{eq:formula for m general field non planar double}
\m(x)=\dfrac{\(Tx-\p(x),g(Tx)-f(\p(x))\)}{\sqrt{|Tx-\p(x)|^2+\(g(Tx)-f(\p(x))\)^2}}.
\end{equation} 
Repeating the steps of Section~\ref{sec:single metasurface} with this $\m$, we obtain
\begin{align}\label{eq:main formula for T2 non flat double}
&\dfrac{Tx-\p(x)+(g(Tx)-f(\p(x)))\nabla f(\p(x))}{\sqrt{(g(Tx)-f(\p(x)))^2+|Tx-\p(x)|^2}}\notag\\
&=\dfrac{-1}{n_2}\(Id+\nabla f\otimes \nabla f\)\(\Phi_{x_1},\Phi_{x_2}\)
+\dfrac{n_1}{n_2}e_3\nabla f+\dfrac{n_1}{n_2}(e_1(x),e_2(x)).
\end{align}

To connect this with optimal transport theory, we introduce the cost function
\begin{equation}\label{eq:cost for double meta}
c(x,y)=\sqrt{(g(y)-f(\p(x)))^2+|\p(x)-y|^2},
\end{equation}
where $(x,y) \in \Omega_0 \times \Omega_1$.
Notice that if $J_\p=\left(\dfrac{\partial \p_i}{\partial x_j}\right)_{i,j=1}^2$ 
denotes the Jacobian matrix of $\p$, then from \eqref{eq:cost for double meta} we have
\begin{equation}\label{eq:formula for Dx of cost double metalens}
\nabla_xc(x,y)=J_\p(x)\,\dfrac{\p(x)-y-(g(y)-f(\p(x)))\nabla f(\p(x))}{\sqrt{(g(y)-f(\p(x)))^2+|\p(x)-y|^2}}.
\end{equation}
Since the map $\p$ is smooth and invertible, we obtain from \eqref{eq:main formula for T2 non flat double} and \eqref{eq:formula for Dx of cost double metalens}
that 
\begin{equation}\label{eq:formula for Dx of cost optimal map double}
\(J_\p(x)\)^{-1}\nabla_xc(x,Tx)=\dfrac{1}{n_2}\(Id+\nabla f\otimes \nabla f\)\(\Phi_{x_1},\Phi_{x_2}\)
-\dfrac{n_1}{n_2}e_3(x)\nabla f - \dfrac{n_1}{n_2}(e_1(x),e_2(x)),
\end{equation}
showing that $T$ and $\Phi$ are related via the cost $c$. As discussed in Section~\ref{sec:single metasurface}, if $T$ is the optimal map for the cost $c$ and the densities $\rho_0,\rho_1$, and if the assumptions of Proposition~\ref{prop:representation formula from injectivity} hold, then the phase $\Phi$ is determined by
\begin{equation}\label{eq:formula for Dx of cost optimal map double BIS}
\(J_\p(x)\)^{-1}\nabla\psi(x)=\dfrac{1}{n_2}\(Id+\nabla f\otimes \nabla f\)\(\Phi_{x_1},\Phi_{x_2}\)
-\dfrac{n_1}{n_2}e_3(x)\nabla f - \dfrac{n_1}{n_2}(e_1(x),e_2(x)),
\end{equation}
where $\nabla_xc(x,Tx) = \nabla\psi(x)$ for a.e. $x\in\Omega_0$ and for some $c$-concave function $\psi$.

Applying the Sherman--Morrison formula as before, we obtain the phase $\Phi$ as  
\begin{equation}\label{eq:formula for the first phase Phi double}
\(\Phi_{x_1},\Phi_{x_2}\)
=
\(Id-\dfrac{\nabla f\otimes \nabla f}{1+|\nabla f|^2} \)\(n_2\,\(J_\p(x)\)^{-1}\nabla_xc(x, Tx)+n_1 e_3(x)\, \nabla f + n_1\,(e_1(x),e_2(x))\),
\end{equation}
with $\Phi_{x_i}$ evaluated at $\(\p(x),f(\p(x))\)$ and $\nabla f$ evaluated at $\p(x)$.

To determine $\Psi$, we apply the GSL on the surface $S_2$:  
\begin{equation*}
    n_2\textbf{m}(x) - n_3\,(0,0,1) = \lambda\nu_g(Tx) + \nabla\Psi(Tx,g(Tx)),
\end{equation*}
where $\textbf{m}(x)=(m_1(x),m_2(x),m_3(x))$ is now the incident unit vector given in 
\eqref{eq:formula for m general field non planar double}, $\nu_g(Tx) = (-\nabla g(Tx),1)$ is the normal vector to $S_2$, and $\lambda\in \R$; here $x\in \Omega_0$. Since $\Psi$ is required to be tangential,
\begin{equation}\label{eq:tangential condition5}
-\Psi_{x_1}(Tx,g(Tx))\,g_{x_1}(Tx)-\Psi_{x_2}(Tx,g(Tx))\,g_{x_2}(Tx)+\Psi_{x_3}(Tx,g(Tx))=0.
\end{equation}
We then have the equations
\begin{equation}\label{eq:equations for lambda and Phi_xi4}
\begin{cases}
n_2\,m_1&=-\lambda\,g_{x_1}+\Psi_{x_1}\\
n_2\,m_2&=-\lambda\,g_{x_2}+\Psi_{x_2}\\
n_2\,m_3-n_3&=\lambda+\Psi_{x_3},
\end{cases}
\end{equation}
with $\nabla g$ evaluated at $Tx$ and $\Psi_{x_i}$ evaluated at $\(Tx,g(Tx)\)$.
Recall that the light ray first travels from $(x,0)$ to $\(\p(x),f(\p(x))\)$ on $S_1$ with unit direction $\be(x)$, then travels from this point to $(Tx, g(Tx))$ on $S_2$, and finally travels to $(Tx,\beta)$ as in Figure \ref{fig:two metasurfaces}. From \eqref{eq:equations for lambda and Phi_xi4} and the form of $\m$ in \eqref{eq:formula for m general field non planar double}, it follows that 
\[
Tx - \p(x) = \sqrt{(g(Tx)-f(\p(x)))^2+|Tx-\p(x)|^2}\dfrac{1}{n_2}\(-\lambda g_{x_1} + \Psi_{x_1}, -\lambda g_{x_2} + \Psi_{x_2}\).
\]
From the last equation in \eqref{eq:equations for lambda and Phi_xi4}, the value of $m_3$, and \eqref{eq:tangential condition5}, we have 
\begin{align*}
\lambda&=-n_3 + n_2\,\dfrac{g(Tx)-f(\p(x))}{\sqrt{(g(Tx)-f(\p(x)))^2+|Tx-\p(x)|^2}}-\Psi_{x_3}\\
&=-n_3+n_2\,\dfrac{g(Tx)-f(\p(x))}{\sqrt{(g(Tx)-f(\p(x)))^2+|Tx-\p(x)|^2}}-\Psi_{x_1}\,g_{x_1}-\Psi_{x_2}\,g_{x_2}.
\end{align*}
Hence
\[
-\lambda\,g_{x_i}+\Psi_{x_i}=
\(n_3-n_2\,\dfrac{g(Tx)-f(\p(x))}{\sqrt{(g(Tx)-f(\p(x)))^2+|Tx-\p(x)|^2}}+\Psi_{x_1}\,g_{x_1}+\Psi_{x_2}\,g_{x_2}\)\,g_{x_i}+\Psi_{x_i}
\]
for $i=1,2$, and so 
\begin{align*}
& -\lambda\,\nabla g +\(\Psi_{x_1},\Psi_{x_2}\) \\
&= 
\(n_3-n_2\,\dfrac{g(Tx)-f(\p(x))}{\sqrt{(g(Tx)-f(\p(x)))^2+|Tx-\p(x)|^2}}\)\nabla g
+
\(\nabla g\otimes \nabla g\)\(\Psi_{x_1},\Psi_{x_2}\)+\(\Psi_{x_1},\Psi_{x_2}\) \\
&=
n_3\,\nabla g
-n_2\,\dfrac{g(Tx)-f(\p(x))}{\sqrt{(g(Tx)-f(\p(x)))^2+|Tx-\p(x)|^2}}\,\nabla g
+
\(Id+\nabla g\otimes \nabla g\)\(\Psi_{x_1},\Psi_{x_2}\).
\end{align*}

Combining these identities, we obtain
\begin{equation}\label{eq:main formula for T2}
\dfrac{Tx-\p(x)+(g(Tx)-f(\p(x)))\nabla g(Tx)}{\sqrt{(g(Tx)-f(\p(x)))^2+|Tx-\p(x)|^2}}
=\dfrac{1}{n_2}\(Id+\nabla g\otimes \nabla g\)\(\Psi_{x_1},\Psi_{x_2}\)
+\dfrac{n_3}{n_2}\nabla g,
\end{equation}
where $\Psi_{x_i}$ is evaluated at $\(Tx,g(Tx)\)$ and $\nabla g$ is evaluated at $Tx$.
From the form of the cost $c$ in \eqref{eq:cost for double meta}, the left-hand side of \eqref{eq:main formula for T2} is $\nabla_yc(x,Tx)$. Using the Sherman--Morrison formula as in Section~\ref{sec:single metasurface}, we obtain
\begin{equation}\label{final Psi}
\(\Psi_{x_1},\Psi_{x_2}\)
=
\(Id-\dfrac{\nabla g\otimes \nabla g}{1+|\nabla g|^2} \)\(n_2\,\nabla_yc(x,Tx)- n_3\,\nabla g\),
\end{equation}
showing that $T$ and $\Psi$ are related via the cost $c$. As before, at this point we do not yet know whether the map $T$ satisfies the conservation condition \eqref{eq:conservation of energy}.
Using the same argument as above, if $T$ is the optimal map for the cost $c$, then $T$ satisfies the conservation condition \eqref{eq:conservation of energy}, and the phase $\Psi$ is given by \eqref{final Psi}.

Let us now examine the relationship between the phases $\Phi$ and $\Psi$.
From \eqref{eq:formula for Dx of cost double metalens}
\[
\(J_\p(x)\)^{-1}\,\nabla_xc(x,Tx)=-\dfrac{Tx-\p(x)+(g(Tx)-f(\p(x)))\nabla f(\p(x))}{\sqrt{(g(Tx)-f(\p(x)))^2+|Tx-\p(x)|^2}}.
\]
Then
\begin{equation}\label{combining x-y derivative}
\begin{split}
\nabla_yc(x,Tx)&=\dfrac{Tx-\p(x)+(g(Tx)-f(\p(x)))\nabla g(Tx)}{\sqrt{(g(Tx)-f(\p(x)))^2+|Tx-\p(x)|^2}}\\
&=
\dfrac{Tx-\p(x)+(g(Tx)-f(\p(x)))\(\nabla g(Tx)-\nabla f(\p(x))+\nabla f(\p(x))\)}{\sqrt{(g(Tx)-f(\p(x)))^2+|Tx-\p(x)|^2}}\\
&=
\dfrac{Tx-\p(x)+(g(Tx)-f(\p(x)))\nabla f(\p(x))}{\sqrt{(g(Tx)-f(\p(x)))^2+|Tx-\p(x)|^2}}
+
\dfrac{(g(Tx)-f(\p(x)))\(\nabla g(Tx)-\nabla f(\p(x))\)}{\sqrt{(g(Tx)-f(\p(x)))^2+|Tx-\p(x)|^2}}\\
&=
-\(J_\p(x)\)^{-1}\,\nabla_xc(x,Tx)
+
\dfrac{(g(Tx)-f(\p(x)))\(\nabla g(Tx)-\nabla f(\p(x))\)}{\sqrt{(g(Tx)-f(\p(x)))^2+|Tx-\p(x)|^2}}.
\end{split}
\end{equation}
Therefore, from \eqref{eq:formula for Dx of cost optimal map double}, the phases $\Phi$ and $\Psi$ satisfy
\begin{align}\label{Psi equation}
\begin{split}
&\dfrac{1}{n_2}\(Id+\nabla g\otimes \nabla g\)\(\Psi_{x_1},\Psi_{x_2}\)
+\dfrac{n_3}{n_2}\nabla g\\
&=
-\dfrac{1}{n_2}\(Id+\nabla f\otimes \nabla f\)\(\Phi_{x_1},\Phi_{x_2}\)
-\dfrac{n_1}{n_2}e_3(x)\nabla f - \dfrac{n_1}{n_2}(e_1(x),e_2(x))\\
&\qquad +
\dfrac{(g(Tx)-f(\p(x)))\(\nabla g(Tx)-\nabla f(\p(x))\)}{\sqrt{(g(Tx)-f(\p(x)))^2+|Tx-\p(x)|^2}},
\end{split}
\end{align}
where $T$ is the optimal map with respect to the cost \eqref{eq:cost for double meta}.

In the particularly important case when $f$ and $g$ are both constant and $\be(x)=(0,0,1)$, we obtain from \eqref{Psi equation} that
\begin{equation}\label{eq:g and f are constant}
\(\Psi_{x_1},\Psi_{x_2}\)
=
- \(\Phi_{x_1},\Phi_{x_2}\)
 ,
\end{equation}
where $\Phi_{x_i}$ is evaluated at $(\varphi(x), f(\varphi(x)))$ and $\Psi_{x_i}$ at $(Tx, g(Tx))$.
%

Additionally, if the assumptions of Proposition~\ref{prop:representation formula from injectivity} hold, then 
from \eqref{final Psi} and \eqref{combining x-y derivative} we obtain
\begin{align}
\begin{split}
&\(\Psi_{x_1},\Psi_{x_2}\) =\\
&\(Id - \dfrac{\nabla g \otimes \nabla g}{1 + |\nabla g|^2}\)\(-n_2\(J_\p(x)\)^{-1}\nabla \psi(x)
+
n_2\dfrac{(g(Tx)-f(\p(x)))\(\nabla g(Tx)-\nabla f(\p(x))\)}{\sqrt{(g(Tx)-f(\p(x)))^2+|Tx-\p(x)|^2}} - n_3\nabla g\),
\end{split}
\end{align}
where $\nabla_xc(x,Tx) = \nabla\psi(x)$ for a.e. $x\in\Omega_0$ and for some $c$-concave function $\psi$.

We have proved the following theorem.

\begin{theorem}\label{thm:main theorem for two metasurfaces}
Let $f, g, \varphi$, and $\textbf{e}(x) = (e_1(x),e_2(x),e_3(x))$ be as in the setup.
If the cost $c$ given by \eqref{eq:cost for double meta} satisfies the assumptions of Proposition~\ref{prop:representation formula from injectivity}\footnote{If $f$ and $g$ satisfy conditions \eqref{twist1}-\eqref{twist2}, then the map $y\mapsto \nabla_xc(x,y)$ is injective for all $x$; see Section~\ref{subsect:twist condition for general g}.}, and \eqref{eq:conservation of energy global} holds with  
$\rho_0$ strictly positive a.e.,
then there is a $c$-concave function $\psi$ 
such that the phases $\Phi$ and $\Psi$ solving the problem described in this section satisfy
\[
\(\Phi_{x_1},\Phi_{x_2}\)
=
\(Id - \dfrac{\nabla f \otimes \nabla f}{1 + |\nabla f|^2}\)\(n_2\,\(J_\p(x)\)^{-1}\nabla\psi(x)+n_1 e_3(x)\, \nabla f + n_1\,(e_1(x),e_2(x))\),
\] 
and 
\begin{equation*}
\begin{split}
&\(\Psi_{x_1},\Psi_{x_2}\) =\\
&\(Id - \dfrac{\nabla g \otimes \nabla g}{1 + |\nabla g|^2}\)\(-n_2\(J_\p(x)\)^{-1}\nabla \psi(x)
+
n_2\dfrac{(g(Tx)-f(\p(x)))\(\nabla g(Tx)-\nabla f(\p(x))\)}{\sqrt{(g(Tx)-f(\p(x)))^2+|Tx-\p(x)|^2}} - n_3\nabla g\),
\end{split}
\end{equation*}
where $\psi$ is the function in \eqref{eq:representation formula for gradient of c}. Here, $\Phi_{x_i}$ is evaluated at $(\varphi(x),f(\varphi(x)))$ and $\nabla f$ is evaluated at $\varphi(x)$; $\Psi_{x_i}$ is evaluated at $(Tx, g(Tx))$ and $\nabla g$ is evaluated at $Tx$.

\end{theorem}

\subsection{Point source case and collimated incident field}
In this subsection, we mention two examples of incident fields that are important and used in optical applications.

The first example is the point-source case. Let $P=(p_1,p_2,p_3)$ be a point below the plane $z=0$, that is, $p_3<0$, and suppose rays emanate from $P$. If $x\in \Omega_0$, then the ray from $P$ to $(x,0)$ has direction $(x_1-p_1,x_2-p_2,-p_3)$, and the field of unit directions at each point $(x,0)\in \Omega_0\times \{0\}$ is
\[
\be(x)=\dfrac{(x_1-p_1,x_2-p_2,-p_3)}{\sqrt{(x_1-p_1)^2+(x_2-p_2)^2+p_3^2}}.
\]
In this case, $\(e_1(x),e_2(x)\)=\nabla_x\(|(x,0)-P|\)$, and the corresponding phase $\Phi$ is given by \eqref{eq:formula for the first phase Phi}.

The second example is the collimated case, where the incident field is $\be(x)=(0,0,1)$.

\setcounter{equation}{0}

\section{Regularity of optimal maps for the cost \texorpdfstring{$\sqrt{(g(y)-f(\p(x)))^2+|\p(x)-y|^2}$}{Lg}}\label{sec:regularity of the optimal map and costs}

In this section, we determine conditions on $f$, $g$, and $\varphi$ so that {\it the partial regularity theory} of De Philippis and Figalli is applicable to the optimal maps appearing in our problems, specifically \cite[Theorem 1.3]{partial-regularity_optimalmaps}, and therefore the optimal maps  appearing in Sections \ref{sec:single metasurface} and \ref{doublemetasurface}
are sufficiently regular on $\Omega$ except possibly on a relatively closed set $\Sigma \subset \Omega$ that has measure zero. 

We begin by recalling the partial regularity theorem of De Philippis and Figalli that will be used below. 

\begin{theorem}[{\cite[Theorem 1.3]{partial-regularity_optimalmaps}\footnote{For a variational proof of this result see \cite[Corollary 1.4]{otto-prodhomme-ried:variationalregularityOTgeneralcosts}.}}]\label{thm:DPF partial regularity}
Let $X,Y\subset \mathbb{R}^n$ be two bounded open sets, and let $\rho_X:X\to \mathbb{R}_+$ and $\rho_Y:Y\to \mathbb{R}_+$ be two continuous probability densities, respectively bounded away from zero and infinity on $X$ and $Y$.
Assume that the cost $c:X\times Y\to \mathbb{R}$ satisfies (C0)--(C3), and denote by $T:X\to Y$ the unique optimal transport map sending $\rho_X$ onto $\rho_Y$.
Then there exist two relatively closed sets $\Sigma_X\subset X$, $\Sigma_Y\subset Y$ of measure zero such that $T:X\setminus \Sigma_X\to Y\setminus \Sigma_Y$ is a homeomorphism of class $C^{0,\beta}_{\mathrm{loc}}$ for any $\beta<1$.

In addition, if $c\in C^{k+2,\alpha}_{\mathrm{loc}}(X\times Y)$, $\rho_X\in C^{k,\alpha}_{\mathrm{loc}}(X)$, and $\rho_Y\in C^{k,\alpha}_{\mathrm{loc}}(Y)$ for some $k\geq 0$ and $\alpha\in(0,1)$, then $T:X\setminus \Sigma_X\to Y\setminus \Sigma_Y$ is a diffeomorphism of class $C^{k+1,\alpha}_{\mathrm{loc}}$.
\end{theorem}

We apply this theorem with $X=\Omega_0$, $Y=\Omega_1$, $\rho_X=\rho_0$ and $\rho_Y=\rho_1$. Therefore, if 
$c\in C^{k+2,\alpha}_{\mathrm{loc}}(\Omega_0\times \Omega_1)$, the densities $\rho_0$, $\rho_1$ are bounded away from zero and infinity, and $\rho_i\in C^{k,\alpha}_{\mathrm{loc}}(\Omega_i)$, $i=0,1$ for some $k\geq 0$, $0<\alpha<1$, then to obtain that the optimal map $T$ is an homeomorphism of class $C^{k+1,\alpha}_{\mathrm{loc}}$, it remains to verify the structural hypotheses (C0)--(C3) for the cost arising in our problems. That cost is
\begin{equation}\label{eq:general formula for the cost}
c(x,y)=\sqrt{(g(y)-f(\p(x)))^2+|\p(x)-y|^2},
\end{equation}
for $(x,y)\in \Omega_0\times\Omega_1$.
Notice that if $f,g,\varphi$ are $C^{k,\alpha}_{\mathrm{loc}}$, $g>f$, it follows from the form of the cost that $c\in C^{k+2,\alpha}_{\mathrm{loc}}(\Omega_0\times \Omega_1)$.

The purpose of this section is to find assumptions on the functions $f$, $g$, and the mapping $\p$ which ensure that this cost satisfies the following conditions:
\begin{enumerate}
   \item[(C0)] $c:\Omega_0 \times \Omega_1 \longrightarrow \mathbb{R}$ is $C^2$ with $\|c\|_{C^2(\Omega_0\times\Omega_1)} < \infty$.
   \item[(C1)] For every $x \in \Omega_0$, the map $\Omega_1 \ni y \longmapsto -D_xc(x,y)$ is injective.
   \item[(C2)] For every $y \in \Omega_1$, the map $\Omega_0 \ni x \longmapsto -D_yc(x,y)$ is injective.
   \item[(C3)] $\det(D_{xy}c)(x,y)\neq 0$ for all $(x,y)\in\Omega_0\times\Omega_1$.
\end{enumerate}

    We remark that when $f$ and $g$ are constant functions, $g>f$, and the incident field is ${\bf e}(x)=(0,0,1)$, the cost $c$  
 is analyzed in \cite[Section 6]{MaTrudingerWang:regularityofpotentials}. There it is proved that $c$ satisfies the (A3)-condition, which, together with $c$-convexity assumptions on the domains $\Omega_0,\Omega_1$ and smoothness of the densities, yields that $T$ is smooth everywhere in $\Omega_0$ \cite[Theorem 2.1]{MaTrudingerWang:regularityofpotentials}; see also \cite{2009trudingerwang:differentiabilitypotentialsOT}.
When $f$ and $g$ are nonconstant, the verification of the (A3)-condition, as well as the requirement of $c$-convexity of the domains, becomes very complicated.

    \subsection{Condition C0} Clearly, $c(x,y)$ is $C^2$ if both $f,g$ and $\p$ are $C^2$ since $g>f$.
    \subsection{Conditions C1 and C2 (Twist Condition)} We seek conditions on the functions $f$ and $g$ such that the twist condition holds.
%
%
 Notice that if 
 \[
 c'(x,y)=\sqrt{(f(x)-g(y))^2+|x-y|^2},
 \] 
 then $c(x,y)=c'(\p(x),y)$. Recall that $\p: \Omega_0 \to \Omega_0'$ is a $C^2$ one-to-one mapping onto another planar domain $\Omega_0'$, and $f \in C^1\(\Omega_0'\)$. Then $\nabla_xc(x,y)=J_{\p}(x)^t\nabla_xc'(\p(x),y)$ and $\nabla_yc(x,y)=\nabla_yc'(\p(x),y)$. Since $\p$ is injective and $J_{\p}(x)$ is nonsingular, it is enough to show that $c'(x, y)$ satisfies the twist condition. To simplify the notation, from now on we write $c(x,y)=\sqrt{(f(x)-g(y))^2+|x-y|^2}$, a cost now defined for $x\in \Omega_0'$ and $y\in \Omega_1$. 
 
 We next analyze the validity of the twist condition separately in two cases: first when $g$ is constant, and second when $g$ is general. In fact, when $g$ is constant, the conditions on $f$ are slightly more general than the corresponding conditions when $g$ is general.
 
 
   
\subsubsection{Verification of the twist condition when $g=\beta$}\label{subsect:verification of twist g=beta}
We aim to show that the cost function 
\[c(x,y) = \sqrt{(\beta-f(x))^2 + |x-y|^2}\]
satisfies the C1 condition under certain assumptions on the function $f(x)$. 
Note that \[\nabla_x c(x,y) = \frac{-\nabla f(x)(\beta-f(x)) + (x-y)}{\sqrt{(\beta-f(x))^2 + |x-y|^2}}.\]
Suppose, for a contradiction, that $\nabla_x c(x,y_1) = \nabla_x c(x,y_2)$ for some $y_1 \neq y_2$, and set 
\[\nabla_x c(x,y_1) = \nabla_x c(x,y_2) := v.\]
Hence
\begin{align*}
-\nabla f(x)(\beta-f(x)) + (x-y_1) &= v\sqrt{(\beta-f(x))^2 + |x-y_1|^2} \\
-\nabla f(x)(\beta-f(x)) + (x-y_2) &= v\sqrt{(\beta-f(x))^2 + |x-y_2|^2}.
\end{align*}
Subtracting these equations yields
\[y_2 - y_1 = v\left(\sqrt{(\beta-f(x))^2 + |x-y_1|^2} - \sqrt{(\beta-f(x))^2 + |x-y_2|^2}\right).
\]
We claim that
\begin{equation}\label{eq:inequality for square roots}
\left|\sqrt{(\beta-f(x))^2 + |x-y_1|^2} - \sqrt{(\beta-f(x))^2 + |x-y_2|^2}\right| < |y_2 - y_1|.
\end{equation}
In fact, let $A = (\beta-f(x))^2$, $B_1 = |x-y_1|^2$, and $B_2 = |x-y_2|^2$. Then
\begin{align*}
\left|\sqrt{A + B_1} - \sqrt{A + B_2}\right| = \frac{|B_1 - B_2|}{\sqrt{A + B_1} + \sqrt{A + B_2}} 
< \frac{|B_1 - B_2|}{\sqrt{B_1} + \sqrt{B_2}}
&\leq |y_2 - y_1|.
\end{align*}
We then get $|y_2 - y_1|<|v|\,|y_2 - y_1|$, and since $y_2\neq y_1$, it follows that $|v|>1$.
We shall then find conditions on $f$ such that $|\nabla_x c(x,y)|\leq 1$ for all $x,y$ which 
yields, as desired, that $\nabla_x c(x,y)$ is injective as a function of $y$. 

The condition $|\nabla_x c(x,y)| \leq1$ is equivalent to 
\[\left|-\nabla f(x)(\beta-f(x)) + (x-y)\right|^2 \leq (\beta-f(x))^2 + |x-y|^2\]
for all $x$ and $y$,
that is,
\[|\nabla f(x)|^2(\beta-f(x))^2  - 2(\beta-f(x))\nabla f(x)\cdot (x-y) \leq (\beta-f(x))^2, 
\]
and since $\beta >f(x)$ we get 
\begin{equation}\label{C1conditionbeta}
|\nabla f(x)|^2(\beta-f(x)) - 2\nabla f(x)\cdot (x-y) \leq  \beta-f(x).
\end{equation}
If we let 
\begin{align*}
G=\max_{x\in \Omega_0',y\in \Omega_1}|x-y|;&\quad
M_0=\min_{x\in \Omega_0',y\in \Omega_1}|f(x)-\beta|>0;\quad
M_f=\sup_{x\in \Omega_0'}|\nabla f(x)|
\end{align*}
and we assume that 
\begin{equation}\label{eq:C1condition for beta in terms of max and min}
M_f^2\(\beta+\|f\|_{L^\infty}\)+2 M_f\, G\leq M_0
\end{equation}
then \eqref{C1conditionbeta} follows, and therefore $\nabla_x c(x,y)$ is injective as a function of $y$.
%


It remains to see when $\nabla_y c(x,y)$ is injective as a function of $x$. We will show that this holds for each $y$ since $M_0>0$.
Note that \[\nabla_y c(x,y) = \frac{y-x}{\sqrt{(\beta-f(x))^2 + |x-y|^2}}.\]
Suppose, for a contradiction, that $\nabla_y c(x_1,y) = \nabla_y c(x_2,y)$ for some $x_1 \neq x_2$, that is, 
\[\frac{y-x_1}{\sqrt{(\beta-f(x_1))^2 + |x_1-y|^2}} = \frac{y-x_2}{\sqrt{(\beta-f(x_2))^2 + |x_2-y|^2}} := v.
\]
Hence
\begin{equation*}
y-x_1 = v\sqrt{(\beta-f(x_1))^2 + |x_1-y|^2},\qquad
y-x_2 = v\sqrt{(\beta-f(x_2))^2 + |x_2-y|^2},
\end{equation*}
and subtracting these equations:
\[x_2 - x_1 = v\left(\sqrt{(\beta-f(x_1))^2 + |x_1-y|^2} - \sqrt{(\beta-f(x_2))^2 + |x_2-y|^2}\right).\]
Set $c_i=\sqrt{(\beta-f(x_i))^2 + |x_i-y|^2}$ and $a_i=\beta-f(x_i)$ for $i=1,2$.
Since $a_i>0$ and $y-x_i=vc_i$, we have $c_i=a_i/\sqrt{1-|v|^2}$. Hence
\[
|x_2-x_1|=\dfrac{|v|}{\sqrt{1-|v|^2}}\,|f(x_2)-f(x_1)|.
\]
Moreover, $|v|/\sqrt{1-|v|^2}=|x_i-y|/(\beta-f(x_i))\leq G/M_0$. Therefore
\[
|x_2-x_1|\leq \dfrac{GM_f}{M_0}|x_2-x_1|.
\]
Since \eqref{eq:C1condition for beta in terms of max and min} implies $GM_f/M_0<1$, this is a contradiction. Hence $\nabla_y c(x,y)$ is injective as a function of $x$.

\subsubsection{Verification of the twist condition for general $g$}\label{subsect:twist condition for general g}
We will find conditions on $f$ and $g$ so that the twist condition holds. Recall that the cost is 
$c(x,y)=\sqrt{(f(x)-g(y))^2+|x-y|^2}$.
We start with the injectivity of $\nabla_x c(x,y)$ as a function of $y$.
We have
\[\nabla_x c(x,y) = \frac{\nabla f(x)(f(x)-g(y)) + x-y}{\sqrt{(f(x)-g(y))^2 + |x-y|^2}}.\]
Assume for a contradiction that $\nabla_x c(x,y_1) = \nabla_x c(x,y_2)$ for some $y_1 \neq y_2$, and let  
\[
\frac{\nabla f(x)(f(x)-g(y_1)) + x-y_1}{\sqrt{(f(x)-g(y_1))^2 + |x-y_1|^2}}
=
\frac{\nabla f(x)(f(x)-g(y_2)) + x-y_2}{\sqrt{(f(x)-g(y_2))^2 + |x-y_2|^2}}:=v.
\]
Then
\begin{align*}
\nabla f(x)(f(x)-g(y_1)) + x-y_1
&=v\,\sqrt{(f(x)-g(y_1))^2 + |x-y_1|^2}\\
\nabla f(x)(f(x)-g(y_2)) + x-y_2
&=v\,\sqrt{(f(x)-g(y_2))^2 + |x-y_2|^2}.
\end{align*}
Subtracting the second from the first equation yields
\[
y_2-y_1+\nabla f(x)\left(g(y_2)-g(y_1)\right)
=
v\left(\sqrt{(f(x)-g(y_1))^2 + |x-y_1|^2}
-
\sqrt{(f(x)-g(y_2))^2 + |x-y_2|^2} \right).
\]
We claim that
\begin{equation}\label{key inequality general cost}
\left|\sqrt{(g(y_1)-f(x))^2 + |x-y_1|^2} - \sqrt{(g(y_2)-f(x))^2 + |x-y_2|^2}\right| < |g(y_2) - g(y_1)| + |y_2 - y_1|.
\end{equation}

To prove \eqref{key inequality general cost}, we use the inequality
\begin{equation}\label{eq:triangle inequality with square root}
\left|\sqrt{A_1 + B_1} - \sqrt{A_2 + B_2}\right|
<
\left|\sqrt{A_1 } - \sqrt{A_2 }\right|
+
\left|\sqrt{ B_1} - \sqrt{ B_2}\right|
\end{equation}
for $A_i>0, B_i\geq 0$.
If we apply this inequality with
\[
A_1 = (f(x)-g(y_1))^2, \quad A_2 = (f(x)-g(y_2))^2, \quad B_1 = |x-y_1|^2,\quad B_2 = |x-y_2|^2, 
\]
we obtain \eqref{key inequality general cost}.
Therefore,
\begin{align*}
&\left|y_2-y_1+\nabla f(x)\left(g(y_2)-g(y_1)\right) \right|\\
&=
|v|\,\left|\sqrt{(f(x)-g(y_1))^2 + |x-y_1|^2}
-
\sqrt{(f(x)-g(y_2))^2 + |x-y_2|^2} \right|\\
&<
|v|\,\left(|g(y_2)-g(y_1)|+|y_2 - y_1|\right).
\end{align*}

For $x\in \Omega_0'$ and $y\in \Omega_1$, let 
\[
v(x,y)=\frac{\nabla f(x)(f(x)-g(y)) + x-y}{\sqrt{(f(x)-g(y))^2 + |x-y|^2}}.
\]
Let 
\begin{align*}
G=\max_{x\in \Omega_0',y\in \Omega_1}|x-y|;&\quad
M_0=\min_{x\in \Omega_0',y\in \Omega_1}|f(x)-g(y)|>0;\\
M_f=\sup_{x\in \Omega_0'}|\nabla f(x)|;
&\quad M_g=\sup_{y\in \Omega_1}|\nabla g(y)|.
\end{align*}
Pick $0<\alpha<1$. If $M_f+\dfrac{G}{M_0}\leq \alpha$,
then
$|v(x,y)|\leq \alpha$.
In fact, 
\begin{align*}
|v(x,y)|
&\leq \frac{|\nabla f(x)|\,|f(x)-g(y)| + |x-y|}{\sqrt{(f(x)-g(y))^2 + |x-y|^2}}\\
&\leq
\frac{|\nabla f(x)|\,|f(x)-g(y)|}{\sqrt{(f(x)-g(y))^2 + |x-y|^2}}
+
\frac{ |x-y|}{\sqrt{(f(x)-g(y))^2 + |x-y|^2}}\\
&\leq
|\nabla f(x)|+\dfrac{|x-y|}{|f(x)-g(y)|}
\leq
M_f+\dfrac{G}{M_0}.
\end{align*}
Hence, we get 
\begin{align*}
\left|y_2-y_1+\nabla f(x)\left(g(y_2)-g(y_1)\right) \right|
&<
\alpha\,\left(|g(y_2)-g(y_1)|+|y_2 - y_1|\right),
\end{align*}
which implies
\begin{align*}
\left||y_2-y_1|-|\nabla f(x)|\left|g(y_2)-g(y_1)\right| \right|
&<
\alpha\,\left(|g(y_2)-g(y_1)|+|y_2 - y_1|\right).
\end{align*}
That is,
\begin{align*}
-\alpha\,\left(|g(y_2)-g(y_1)|+|y_2 - y_1|\right)
&<
|y_2-y_1|-|\nabla f(x)|\left|g(y_2)-g(y_1)\right|\\
&<
\alpha\,\left(|g(y_2)-g(y_1)|+|y_2 - y_1|\right).
\end{align*}
The last inequality reads
\[
(1-\alpha)|y_2-y_1|< \left(|\nabla f(x)|+\alpha\right)|g(y_2)-g(y_1)|
\leq
 \left(M_f+\alpha\right)M_g|y_2-y_1|.
\]
If $y_1\neq y_2$, we get 
\[
(1-\alpha) 
<
 \left(M_f+\alpha\right)M_g,
\]
and so if 
\[
(1-\alpha) 
\geq
 \left(M_f+\alpha\right)M_g,
\]
we get a contradiction.

Therefore, if we fix $0<\alpha<1$ and the functions $f,g$, and the domains $\Omega_0',\Omega_1$ satisfy
\[
M_f+\dfrac{G}{M_0}\leq \alpha,
\quad (1-\alpha) 
\geq
 \left(M_f+\alpha\right)M_g,
\]
then the map $y\mapsto \nabla_xc(x,y)$ is injective in $\Omega_1$ for each fixed $x\in \Omega_0'$.
This means that, assuming smallness conditions on the gradients of $f$ and $g$ and sufficient separation between their graphs, injectivity follows.

The C2 condition, i.e., the injectivity of the map $x\mapsto \nabla_yc(x,y)$ in $\Omega_0'$, follows in the same way. Calculating $\nabla_yc(x,y)$ gives

\[\nabla_y c(x,y) = \frac{\nabla g(y)(g(y)-f(x)) + y - x}{\sqrt{(f(x)-g(y))^2 + |x-y|^2}}.\]
As before, assume for a contradiction that $\nabla_y c(x_1,y) = \nabla_y c(x_2,y)$ for some $x_1 \neq x_2$. 

If 
\[
\frac{\nabla g(y)(g(y)-f(x_1)) + y-x_1}{\sqrt{(f(x_1)-g(y))^2 + |x_1-y|^2}}
=
\frac{\nabla g(y)(g(y)-f(x_2)) + y-x_2}{\sqrt{(f(x_2)-g(y))^2 + |x_2-y|^2}}:=v,
\]
then
\begin{align*}
\nabla g(y)(g(y)-f(x_1)) + y-x_1
&=v\,\sqrt{(f(x_1)-g(y))^2 + |x_1-y|^2}\\
\nabla g(y)(g(y)-f(x_2)) + y-x_2
&=v\,\sqrt{(f(x_2)-g(y))^2 + |x_2-y|^2}.
\end{align*}
Subtracting the second from the first equation yields
\[
x_2-x_1+\nabla g(y)\left(f(x_2)-f(x_1)\right)
=
v\left(\sqrt{(f(x_1)-g(y))^2 + |x_1-y|^2}
-
\sqrt{(f(x_2)-g(y))^2 + |x_2-y|^2} \right).
\]

By \eqref{key inequality general cost},
\[\left|\sqrt{(f(x_1)-g(y))^2 + |x_1-y|^2}
-
\sqrt{(f(x_2)-g(y))^2 + |x_2-y|^2}\right| < |f(x_2)-f(x_1)| + |x_2-x_1|.\]

Hence, 
\begin{align*}
&\left|x_2-x_1+\nabla g(y)\left(f(x_2)-f(x_1)\right) \right|\\
&=
|v|\,\left|\sqrt{(f(x_1)-g(y))^2 + |x_1-y|^2}
-
\sqrt{(f(x_2)-g(y))^2 + |x_2-y|^2} \right|\\
&<
|v|\,\left(|f(x_2)-f(x_1)|+|x_2 - x_1|\right).
\end{align*}

Let $G, M_0, M_f$ and $M_g$ be as before, and \[
v=v(x,y)=\frac{\nabla g(y)(g(y)-f(x)) + y-x}{\sqrt{(f(x)-g(y))^2 + |x-y|^2}},
\]
for $x\in\Omega_0'$ and $y\in\Omega_1$. 

Pick $0<\alpha<1$. If $M_g+\dfrac{G}{M_0}\leq \alpha$,
then
$|v(x,y)|\leq \alpha$.
In fact, 
\begin{align*}
|v(x,y)|
&\leq \frac{|\nabla g(y)|\,|g(y)-f(x)| + |x-y|}{\sqrt{(f(x)-g(y))^2 + |x-y|^2}}\\
&\leq
\frac{|\nabla g(y)|\,|g(y)-f(x)|}{\sqrt{(f(x)-g(y))^2 + |x-y|^2}}
+
\frac{ |x-y|}{\sqrt{(f(x)-g(y))^2 + |x-y|^2}}\\
&\leq
|\nabla g(y)|+\dfrac{|x-y|}{|f(x)-g(y)|}
\leq
M_g+\dfrac{G}{M_0}.
\end{align*}
Hence, we get 
\begin{align*}
\left|x_2-x_1+\nabla g(y)\left(f(x_2)-f(x_1)\right) \right|
&<
\alpha\,\left(|f(x_2)-f(x_1)|+|x_2 - x_1|\right),
\end{align*}
which implies
\begin{align*}
\left||x_2-x_1|-|\nabla g(y)|\left|f(x_2)-f(x_1)\right| \right|
&<
\alpha\,\left(|f(x_2)-f(x_1)|+|x_2 - x_1|\right).
\end{align*}
That is, 
\begin{align*}
-\alpha\,\left(|f(x_2)-f(x_1)|+|x_2 - x_1|\right) &<
|x_2-x_1|-|\nabla g(y)|\left|f(x_2)-f(x_1)\right|\\
&
<
\alpha\,\left(|f(x_2)-f(x_1)|+|x_2 - x_1|\right).
\end{align*}
The last inequality reads
\[
(1-\alpha)|x_2-x_1|< \left(|\nabla g(y)|+\alpha\right)|f(x_2)-f(x_1)|
\leq
 \left(M_g+\alpha\right)M_f|x_2-x_1|.
\]
If $x_1\neq x_2$, we get 
\[
(1-\alpha) 
<
 \left(M_g+\alpha\right)M_f,
\]
and so if 
\[
(1-\alpha) 
\geq
 \left(M_g+\alpha\right)M_f,
\]
we get a contradiction.

Therefore, if we fix $0<\alpha<1$ and the functions $f,g$, and the domains $\Omega_0',\Omega_1$ satisfy
\[
M_g+\dfrac{G}{M_0}\leq \alpha,
\quad (1-\alpha) 
\geq
 \left(M_g+\alpha\right)M_f,
\]
then the map $x\mapsto \nabla_yc(x,y)$ is injective in $\Omega_0'$ for each fixed $y\in \Omega_1$, i.e., the cost $c(x,y) = \sqrt{(f(x)-g(y))^2+|x-y|^2}$ satisfies the C2 condition.

{\it Summarizing}, the cost $c(x,y) = \sqrt{(f(x)-g(y))^2+|x-y|^2}$ satisfies the twist condition if 
\begin{align}
M_f + \dfrac{G}{M_0} &\leq \alpha \label{twist1}, \\
M_g + \dfrac{G}{M_0} &\leq \alpha, \\
(M_f + \alpha)M_g &\leq (1-\alpha), \\
(M_g + \alpha)M_f &\leq (1-\alpha), \label{twist2}
\end{align}
where $0<\alpha<1$ is fixed, and 
 \begin{align*}
G=\max_{x\in \Omega_0',y\in \Omega_1}|x-y|;&\quad
M_0=\min_{x\in \Omega_0',y\in \Omega_1}|f(x)-g(y)|>0;\\
M_f=\sup_{x\in \Omega_0'}|\nabla f(x)|;
&\quad M_g=\sup_{y\in \Omega_1}|\nabla g(y)|.
\end{align*}

Also, when $g=\beta$, the twist condition holds if \eqref{eq:C1condition for beta in terms of max and min} holds.

Qualitatively, the twist condition holds under smallness conditions on the gradients of $f$ and $g$ and sufficient separation between their graphs.

\subsection{Condition C3} 
We seek conditions on the functions $f$, $g$, and the mapping $\p$ so that $\det \dfrac{\partial^2 c}{\partial x_i\partial y_j}\neq 0$. Again, if $c'(x,y)=\sqrt{(f(x)-g(y))^2+|x-y|^2}$, then $c(x,y)=c'(\p(x),y)$. Written in matrix form, the chain rule gives
\[
\dfrac{\partial^2 c}{\partial x\partial y}(x,y)=\dfrac{\partial \p}{\partial x}(x)^t\,\dfrac{\partial^2 c'}{\partial x\partial y}(\p(x),y).
\] 
If $\p$ has nonsingular Jacobian, then $\det \dfrac{\partial \p}{\partial x}(x)\neq 0$.
Therefore, showing that $\det \dfrac{\partial^2 c}{\partial x\partial y}(x,y)\neq 0$ is equivalent to seeking conditions on $f$ and $g$ such that $\det \dfrac{\partial^2 c'}{\partial x\partial y}(\p(x),y)\neq 0$.

As in the previous subsection, we continue using the notation
\[
c(x,y)=\sqrt{(f(x)-g(y))^2+|x-y|^2}.
\]
To simplify the notation, write $\Delta=(f(x)-g(y))^2+|x-y|^2$.
We have
\[
\dfrac{\partial c}{\partial x_i}=\Delta^{-1/2}\((f(x)-g(y))\,\partial_{x_i}f +x_i-y_i\).
\]
So
\begin{equation*}
\resizebox{1\textwidth}{!}{
$\begin{aligned}
\dfrac{\partial^2 c}{\partial x_i\partial y_j}
&=
(-1/2)\Delta^{-3/2}\(-2\,(f(x)-g(y))\,\partial_{y_j}g-2\,(x_j-y_j)\)\,\((f(x)-g(y))\,\partial_{x_i}f +x_i-y_i\)\\
&\qquad + \Delta^{-1/2}\(-\partial_{y_j}g \,\partial_{x_i}f-\delta_{ij}\)\\
&=
\Delta^{-1/2}\(\Delta^{-1}\((f(x)-g(y))\partial_{y_j}g+x_j-y_j\)\,\((f(x)-g(y))\,\partial_{x_i}f +x_i-y_i\)+\(-\partial_{y_j}g \,\partial_{x_i}f-\delta_{ij}\)\),
\end{aligned}$
}
\end{equation*}
that is,
\begin{equation}\label{eq:Jacobian c}
\resizebox{0.91\textwidth}{!}{
$\dfrac{\partial^2 c}{\partial x\partial y}
=
\dfrac{-1}{c(x,y)}
\(\(\dfrac{(f(x)-g(y))\nabla g+x-y}{c(x,y)} \)\otimes \(\dfrac{(g(y)-f(x))\nabla f+y-x}{c(x,y)} \)+\(\nabla g\otimes \nabla f\)+Id \);$
}
\end{equation}
with $u\otimes v=uv^t$.

\subsubsection{Verification of the C3 condition when $g = \beta$}\label{subsect:condition C3 for g=beta}
In this case the surface $S_2$ is a horizontal plane. Then
\begin{equation}\label{eq:hessian of c beta}
\dfrac{\partial^2 c}{\partial x\partial y}
=
\dfrac{-1}{c(x,y)}
\(\(\dfrac{x-y}{c(x,y)} \)\otimes \(\dfrac{(\beta-f(x))\nabla f(x)+y-x}{c(x,y)} \)+Id \).
\end{equation}

Applying the determinant formula \eqref{eq:Sherman Morrinson formula}, it follows that
\begin{align*}
\det \dfrac{\partial^2 c}{\partial x\partial y}
&=\(\dfrac{-1}{c(x,y)}\)^n\,\(1+\dfrac{x-y}{c(x,y)}\cdot \dfrac{(\beta-f(x))\nabla f(x)+y-x}{c(x,y)}\)\\
&=\(\dfrac{-1}{c(x,y)}\)^n\,\dfrac{1}{c(x,y)^2}\(c(x,y)^2+
(\beta-f(x))\,(x-y)\cdot \nabla f(x)-|x-y|^2\)\\
&=
\(\dfrac{-1}{c(x,y)}\)^n\,\dfrac{1}{c(x,y)^2}
\((\beta-f(x))^2+(\beta-f(x))\,(x-y)\cdot \nabla f(x)\)\\
&=
\(\dfrac{-1}{c(x,y)}\)^n\,\dfrac{\beta-f(x)}{c(x,y)^2}
\(\beta-f(x)+(x-y)\cdot \nabla f(x)\).
\end{align*}
Since by assumption the distance between the graph of $f$ and the plane $g=\beta$ is strictly positive, the quantity 
$\(\dfrac{-1}{c(x,y)}\)^n\,\dfrac{\beta-f(x)}{c(x,y)^2}\neq 0$. Thus we need $f$ to satisfy $|\beta-f(x)+(x-y)\cdot \nabla f(x)|\neq 0$. We have
\begin{align*}
|\beta-f(x)+(x-y)\cdot \nabla f(x)|
&\geq \beta-f(x)-|(x-y)\cdot \nabla f(x)|\\
&\geq \beta-f(x)-|x-y|\,|\nabla f(x)|\\
&\geq \beta-f(x)-\max_{x\in \Omega_0',y\in \Omega_1}|x-y|\,\max_{x\in \Omega_0'}|\nabla f(x)|.
\end{align*}
For example, if $f$ satisfies 
\begin{equation}\label{eq:condition C3 on f when g=beta}
\min_{x\in \Omega_0'}\(\beta-f(x)\)>\max_{x\in \Omega_0',y\in \Omega_1}|x-y|\,\max_{x\in \Omega_0'}|\nabla f(x)|,
\end{equation}
then $\det \dfrac{\partial^2 c}{\partial x\partial y}\neq 0$.

\subsubsection{Verification of the C3 condition for general $g$}
For the general case, when $g$ is not necessarily constant, use equation \eqref{eq:Jacobian c}. Let 
\[
A=Id+\nabla g(y)\otimes \nabla f(x)
\]
so $A$ is invertible iff $1+\nabla g(y)\cdot \nabla f(x)\neq 0$, and from \eqref{eq:Sherman Morrinson formula}
\[
A^{-1}=Id -\dfrac{\nabla g(y)\otimes \nabla f(x)}{1+\nabla g(y)\cdot \nabla f(x)}.
\]
Set 
\[
u=\dfrac{(f(x)-g(y))\nabla g(y)+x-y}{c(x,y)},
\qquad 
v=\dfrac{(g(y)-f(x))\nabla f(x)+y-x}{c(x,y)}.
\]
From the determinant formula \eqref{eq:Sherman Morrinson formula} we have
\[
\det \dfrac{\partial^2 c}{\partial x\partial y}
=
\(\dfrac{-1}{c(x,y)}\)^n
\(1+u^t A^{-1}v\)\,\(1+\nabla g(y)\cdot \nabla f(x)\).
\]
We have
\begin{equation*}
\begin{split}
1+u^t A^{-1}v
&=1+u^t\left(Id -\dfrac{\nabla g(y)\otimes \nabla f(x)}{1+\nabla g(y)\cdot \nabla f(x)}\right)v \\
&= 1+u^tv- \dfrac{u^t\nabla g(y) \nabla f(x)^tv}{1+\nabla g(y)\cdot \nabla f(x)} \\
&= 1 + u\cdot v - \dfrac{(u \cdot\nabla g(y)) (\nabla f(x)\cdot v)}{1+\nabla g(y)\cdot \nabla f(x)} \\
&= \dfrac{(1+\nabla g(y)\cdot \nabla f(x))(1+u\cdot v) - (u \cdot\nabla g(y)) (\nabla f(x)\cdot v)}{1+\nabla g(y)\cdot \nabla f(x)}.
\end{split}
\end{equation*}
Also, with $f, \nabla f$ evaluated at $x$ and $g,\nabla g$ evaluated at $y$, we have 
\begin{equation*}  
\resizebox{1\textwidth}{!}{
$\begin{aligned}
u\cdot v
&=
\dfrac{1}{c(x,y)^2}
\left(
-(f-g)^2\nabla g \cdot \nabla f-(f-g)\nabla g\cdot (x-y)-(f-g)(x-y)\cdot\nabla f-(x-y)\cdot(x-y)
\right)
\\
&=
\dfrac{1}{c(x,y)^2}
\left(
-(f-g)^2\nabla g\cdot \nabla f-(f-g)(x-y)\left(\nabla g+\nabla f\right)-|x-y|^2
\right), 
\end{aligned}$
}
\end{equation*}
and
\vspace{-0.5pt}
\begin{align*}
1+u\cdot v
&=
\dfrac{c(x,y)^2-(f-g)^2\nabla g\cdot \nabla f-(f-g)(x-y)\cdot\left(\nabla g+\nabla f\right)-|x-y|^2}{c(x,y)^2}\\
&= \dfrac{(f-g)^2\left(1-\nabla g\cdot \nabla f\right)-(f-g)(x-y)\cdot\left(\nabla g+\nabla f\right)}{c(x,y)^2}.
\end{align*}
Therefore,
\vspace{-0.5pt}
\begin{align*}
&(1+\nabla g(y)\cdot \nabla f(x))(1+u\cdot v) =  \\
&\dfrac{(f-g)^2\left(1-\nabla g\cdot \nabla f\right)(1+\nabla g\cdot \nabla f)-(f-g)(x-y)\cdot\left(\nabla g+\nabla f\right)(1+\nabla g\cdot \nabla f)} {c(x,y)^2} = \\
&\dfrac{(f-g)^2\left(1-(\nabla g\cdot \nabla f)^2\right)-(f-g)(x-y)\cdot\left(\nabla g+\nabla f\right)(1+\nabla g\cdot \nabla f)} {c(x,y)^2} = \\
&\dfrac{(g-f)^2\left(1-(\nabla g\cdot \nabla f)^2\right)+(g-f)(x-y)\cdot\left(\nabla g+\nabla f\right)(1+\nabla g\cdot \nabla f)} {c(x,y)^2},
\end{align*}

and
\begin{equation*}
\resizebox{1\textwidth}{!}{
$\begin{aligned}
(u \cdot\nabla g(y)) (\nabla f(x)\cdot v)
&= \(\ds\frac{(f-g)\nabla g \cdot\nabla g + (x-y)\cdot\nabla g)}{c(x,y)}\)\(\ds\frac{-(f-g)\nabla f\cdot\nabla f - (x-y)\cdot \nabla f}{c(x,y)}\) \\
&= \ds\frac{1}{c(x,y)^2}\((f-g)|\nabla g|^2 + (x-y)\cdot\nabla g\)
\(-(f-g)|\nabla f|^2 - (x-y)\cdot\nabla f\) \\
&= \ds\frac{1}{c(x,y)^2} \bigl(-(f-g)^2|\nabla g|^2|\nabla f|^2 - (f-g)|\nabla g|^2(x-y)\cdot\nabla f \\
& \hspace{2cm} - (f-g)|\nabla f|^2(x-y)\cdot\nabla g - ((x-y)\cdot\nabla g)((x-y)\cdot\nabla f)\bigr) \\
&=\ds\frac{1}{c(x,y)^2} \bigl(-(g-f)^2|\nabla g|^2|\nabla f|^2 + (g-f)|\nabla g|^2(x-y)\cdot\nabla f \\
& \hspace{2cm} + (g-f)|\nabla f|^2(x-y)\cdot\nabla g - ((x-y)\cdot\nabla g)((x-y)\cdot\nabla f)\bigr).
\end{aligned}$
}
\end{equation*}
\newpage
Then
\begin{equation*}
\resizebox{1\textwidth}{!}{
$\begin{aligned}
 &(1+\nabla g(y)\cdot \nabla f(x))(1+u\cdot v) - (u \cdot\nabla g(y)) (\nabla f(x)\cdot v) = \\
 &\dfrac{1}{c(x,y)^2}\bigl((g-f)^2(1-(\nabla g\cdot \nabla f)^2) + (g-f)^2|\nabla g|^2|\nabla f|^2 + (g-f)(x-y)\cdot(\nabla g + \nabla f)(1+\nabla g\cdot\nabla f) \\ 
 & - (g-f)|\nabla g|^2(x-y)\cdot\nabla f - (g-f)|\nabla f|^2(x-y)\cdot\nabla g + ((x-y)\cdot\nabla g)((x-y)\cdot\nabla f)\bigr)  \\
 &= \dfrac{1}{c(x,y)^2}\Bigl((g-f)^2\bigl[1-(\nabla g\cdot \nabla f)^2 + |\nabla g|^2|\nabla f|^2\bigr] + (g-f)[(x-y)\cdot\nabla f][1+\nabla g \cdot \nabla f-|\nabla g|^2] \\
 &+ (g-f)[(x-y)\cdot\nabla g][1+\nabla g \cdot \nabla f-|\nabla f|^2] + ((x-y)\cdot\nabla g)((x-y)\cdot\nabla f)\Bigr) \\
 & \geq \dfrac{1}{c(x,y)^2}\Bigl((g-f)^2 + (g-f)[(x-y)\cdot\nabla f][1+\nabla g \cdot \nabla f-|\nabla g|^2]\\ 
 &+ (g-f)[(x-y)\cdot\nabla g][1+\nabla g \cdot \nabla f-|\nabla f|^2] 
 + ((x-y)\cdot\nabla g)((x-y)\cdot\nabla f)\Bigr):=(*), \text{\small as $(\nabla g\cdot \nabla f)^2 \leq |\nabla g|^2|\nabla f|^2$.} \hspace{0.1cm} 
 \end{aligned}$
 }
\end{equation*}

Let $\alpha_1,\alpha_2$ be positive constants and suppose that 
\begin{equation}\label{eq:condition on f and g giving determinant not zero}
\text{$|(x-y)\cdot\nabla g|\leq \alpha_1\,(g-f)$, $|(x-y)\cdot\nabla f|\leq \alpha_1\,(g-f)$, and $|\nabla f|+|\nabla g|\leq \alpha_2$,}
\end{equation}
where $x\in \Omega_0'$ and $y\in \Omega_1$.
Then, write 
 \[
 (*) := \dfrac{1}{c(x,y)^2}\left((g-f)^2 + A + B + C\right),
 \] 
with
\begin{align*}
    A &= (g-f)\,\((x-y)\cdot\nabla f\)\,\(1+\nabla g \cdot \nabla f-|\nabla g|^2\) \\
    B &= (g-f)\,\((x-y)\cdot\nabla g\)\,\(1+\nabla g \cdot \nabla f-|\nabla f|^2\) \\
    C &= \((x-y)\cdot\nabla g\)\,\((x-y)\cdot\nabla f\).
\end{align*}
From \eqref{eq:condition on f and g giving determinant not zero}, we then have the following inequalities:
\begin{align*}
    A & \geq -\alpha_1\,(g-f)^2\,\(1+|\nabla g||\nabla f|+|\nabla g|^2\)\\
    B & \geq -\alpha_1\,(g-f)^2\,\(1+|\nabla g||\nabla f|+|\nabla f|^2\) \\
    C & \geq -\alpha_1^2 \,(g-f)^2.
\end{align*} 
Hence, 
\begin{align*}
&\dfrac{1}{c(x,y)^2}\left((g-f)^2 + A + B + C\right)\\
&\geq 
\dfrac{1}{c(x,y)^2}
\,(g-f)^2\,
\(1-\alpha_1\,\(1+|\nabla g||\nabla f|+|\nabla g|^2\)-\alpha_1\,\(1+|\nabla g||\nabla f|+|\nabla f|^2\)-\alpha_1^2\)\\
&=
\dfrac{1}{c(x,y)^2}
\,(g-f)^2\,
\(1-2\,\alpha_1\,\(1+|\nabla g||\nabla f|\)-\alpha_1\,\(|\nabla g|^2+|\nabla f|^2\)-
\alpha_1^2\)\\
&=
\dfrac{1}{c(x,y)^2}
\,(g-f)^2\,
\(1-\alpha_1\,\(2+2\,|\nabla g||\nabla f|+|\nabla g|^2+|\nabla f|^2\)-\alpha_1^2\)\\
&=
\dfrac{1}{c(x,y)^2}
\,(g-f)^2\,
\(1-\alpha_1\,\(2+\(|\nabla g|+|\nabla f|\)^2\)-\alpha_1^2\)\\
&\geq 
\dfrac{1}{c(x,y)^2}
\,(g-f)^2\,
\(1-\alpha_1\,\(2+\alpha_2^2\)-\alpha_1^2\)\\
&= 
\dfrac{1}{c(x,y)^2}
\,(g-f)^2\,
\(1-\alpha_1\,\(2+\alpha_2^2+\alpha_1\)\).
\end{align*}
If we assume that $f$ and $g$ satisfy \eqref{eq:condition on f and g giving determinant not zero} with positive $\alpha_1,\alpha_2$ satisfying
\[
1-\alpha_1\,\(2+\alpha_2^2+\alpha_1\)>0,
\]
then $\det \dfrac{\partial^2 c}{\partial x\partial y}\neq 0$.

In summary, the cost $c(x,y) = \sqrt{(f(x)-g(y))^2+|x-y|^2}$ satisfies the C3 condition if $f$ and $g$ satisfy \eqref{eq:condition on f and g giving determinant not zero}, where $\alpha_1,\alpha_2$ are positive and satisfy the preceding inequality. Equivalently, it is enough to assume

\begin{align}
M_f + M_g &\leq \alpha_2 \label{C31}, \\
GM_f &\leq \alpha_1M_0 , \\
GM_g &\leq \alpha_1M_0 \label{C32},
\end{align}
where $\alpha_1, \alpha_2>0$ are such that $1-\alpha_1\,\(2+\alpha_2^2+\alpha_1\)>0$, and 
 \begin{align*}
G=\max_{x\in \Omega_0',y\in \Omega_1}|x-y|;&\quad
M_0=\min_{x\in \Omega_0',y\in \Omega_1}|f(x)-g(y)|>0;\\
M_f=\sup_{x\in \Omega_0'}|\nabla f(x)|;
&\quad M_g=\sup_{y\in \Omega_1}|\nabla g(y)|.
\end{align*}
This means that the C3 condition holds under smallness conditions on the gradients of $f$ and $g$ and sufficient separation between their graphs.

We have proved the following theorem.

\begin{theorem}\label{thm:main result of regularity for general costs}
Let $\Omega_0, \Omega_1 \subset \mathbb{R}^2$ be compact domains, and let $\rho_0:
\Omega_0\to\mathbb{R}_+$ and $\rho_1:\Omega_1\to\mathbb{R}_+$ be continuous, and bounded away from zero and infinity. Let
\[
c(x,y)=\sqrt{(f(\p(x))-g(y))^2+|\p(x)-y|^2},
\qquad (x,y)\in \Omega_0\times\Omega_1,
\]
where $\p:\Omega_0\to\Omega_0'$ is one-to-one with nonsingular Jacobian. Assume that $f$, $g$, and $\p$ are $C^2$, and that $f$ and $g$ satisfy \eqref{twist1}--\eqref{twist2} and \eqref{C31}--\eqref{C32}. Denote by $T:\Omega_0\to\Omega_1$ the unique optimal transport map sending $\rho_0$ onto $\rho_1$ with respect to the cost $c$. Then there exist two relatively closed sets $\Sigma_0\subset\Omega_0$ and $\Sigma_1\subset\Omega_1$ of measure zero such that
\[
T:\Omega_0\setminus\Sigma_0\to\Omega_1\setminus\Sigma_1
\]
is a homeomorphism of class $C^{0,\beta}_{\mathrm{loc}}$ for any $\beta<1$.

In addition, if $f$, $g$, and $\p$ are of class $C^{k+2,\alpha}_{\mathrm{loc}}$, and $\rho_i\in C^{k,\alpha}_{\mathrm{loc}}(\Omega_i)$, $i=0,1$, for some $k\geq 0$ and $\alpha\in(0,1)$, then
\[
T:\Omega_0\setminus\Sigma_0\to\Omega_1\setminus\Sigma_1
\]
is a diffeomorphism of class $C^{k+1,\alpha}_{\mathrm{loc}}$.

\end{theorem}

\newcommand{\etalchar}[1]{$^{#1}$}
\providecommand{\bysame}{\leavevmode\hbox to3em{\hrulefill}\thinspace}
\providecommand{\MR}{\relax\ifhmode\unskip\space\fi MR }
\providecommand{\MRhref}[2]{%
  \href{http://www.ams.org/mathscinet-getitem?mr=#1}{#2}
}
\providecommand{\href}[2]{#2}

\end{document}